\catcode`\@=11

\font\tenmsa=msam10
\font\sevenmsa=msam7
\font\fivemsa=msam5
\font\tenmsb=msbm10
\font\sevenmsb=msbm7
\font\fivemsb=msbm5
\newfam\msafam
\newfam\msbfam
\textfont\msafam=\tenmsa  \scriptfont\msafam=\sevenmsa
  \scriptscriptfont\msafam=\fivemsa
\textfont\msbfam=\tenmsb  \scriptfont\msbfam=\sevenmsb
  \scriptscriptfont\msbfam=\fivemsb

\def\hexnumber@#1{\ifnum#1<10 \number#1\else
 \ifnum#1=10 A\else\ifnum#1=11 B\else\ifnum#1=12 C\else
 \ifnum#1=13 D\else\ifnum#1=14 E\else\ifnum#1=15 F\fi\fi\fi\fi\fi\fi\fi}

\def\msa@{\hexnumber@\msafam}
\def\msb@{\hexnumber@\msbfam}
\mathchardef\boxdot="2\msa@00
\mathchardef\boxplus="2\msa@01
\mathchardef\boxtimes="2\msa@02
\mathchardef\square="0\msa@03
\mathchardef\blacksquare="0\msa@04
\mathchardef\centerdot="2\msa@05
\mathchardef\lozenge="0\msa@06
\mathchardef\blacklozenge="0\msa@07
\mathchardef\circlearrowright="3\msa@08
\mathchardef\circlearrowleft="3\msa@09
\mathchardef\rightleftharpoons="3\msa@0A
\mathchardef\leftrightharpoons="3\msa@0B
\mathchardef\boxminus="2\msa@0C
\mathchardef\Vdash="3\msa@0D
\mathchardef\Vvdash="3\msa@0E
\mathchardef\vDash="3\msa@0F
\mathchardef\twoheadrightarrow="3\msa@10
\mathchardef\twoheadleftarrow="3\msa@11
\mathchardef\leftleftarrows="3\msa@12
\mathchardef\rightrightarrows="3\msa@13
\mathchardef\upuparrows="3\msa@14
\mathchardef\downdownarrows="3\msa@15
\mathchardef\upharpoonright="3\msa@16

\mathchardef\downharpoonright="3\msa@17
\mathchardef\upharpoonleft="3\msa@18
\mathchardef\downharpoonleft="3\msa@19
\mathchardef\rightarrowtail="3\msa@1A
\mathchardef\leftarrowtail="3\msa@1B
\mathchardef\leftrightarrows="3\msa@1C
\mathchardef\rightleftarrows="3\msa@1D
\mathchardef\Lsh="3\msa@1E
\mathchardef\Rsh="3\msa@1F
\mathchardef\rightsquigarrow="3\msa@20
\mathchardef\leftrightsquigarrow="3\msa@21
\mathchardef\looparrowleft="3\msa@22
\mathchardef\looparrowright="3\msa@23
\mathchardef\circeq="3\msa@24
\mathchardef\succsim="3\msa@25
\mathchardef\gtrsim="3\msa@26
\mathchardef\gtrapprox="3\msa@27
\mathchardef\multimap="3\msa@28
\mathchardef\therefore="3\msa@29
\mathchardef\because="3\msa@2A
\mathchardef\doteqdot="3\msa@2B

\mathchardef\triangleq="3\msa@2C
\mathchardef\precsim="3\msa@2D
\mathchardef\lesssim="3\msa@2E
\mathchardef\lessapprox="3\msa@2F
\mathchardef\eqslantless="3\msa@30
\mathchardef\eqslantgtr="3\msa@31
\mathchardef\curlyeqprec="3\msa@32
\mathchardef\curlyeqsucc="3\msa@33
\mathchardef\preccurlyeq="3\msa@34
\mathchardef\leqq="3\msa@35
\mathchardef\leqslant="3\msa@36
\mathchardef\lessgtr="3\msa@37
\mathchardef\backprime="0\msa@38
\mathchardef\risingdotseq="3\msa@3A
\mathchardef\fallingdotseq="3\msa@3B
\mathchardef\succcurlyeq="3\msa@3C
\mathchardef\geqq="3\msa@3D
\mathchardef\geqslant="3\msa@3E
\mathchardef\gtrless="3\msa@3F
\mathchardef\sqsubset="3\msa@40
\mathchardef\sqsupset="3\msa@41
\mathchardef\trianglerighteq="3\msa@44
\mathchardef\trianglelefteq="3\msa@45
\mathchardef\bigstar="0\msa@46
\mathchardef\between="3\msa@47
\mathchardef\blacktriangledown="0\msa@48
\mathchardef\blacktriangleright="3\msa@49
\mathchardef\blacktriangleleft="3\msa@4A
\mathchardef\blacktriangle="0\msa@4E
\mathchardef\triangledown="0\msa@4F
\mathchardef\eqcirc="3\msa@50
\mathchardef\lesseqgtr="3\msa@51
\mathchardef\gtreqless="3\msa@52
\mathchardef\lesseqqgtr="3\msa@53
\mathchardef\gtreqqless="3\msa@54
\mathchardef\Rrightarrow="3\msa@56
\mathchardef\Lleftarrow="3\msa@57
\mathchardef\veebar="2\msa@59
\mathchardef\barwedge="2\msa@5A
\mathchardef\doublebarwedge="2\msa@5B
\mathchardef\angle="0\msa@5C
\mathchardef\measuredangle="0\msa@5D
\mathchardef\sphericalangle="0\msa@5E
\mathchardef\varpropto="3\msa@5F
\mathchardef\smallsmile="3\msa@60
\mathchardef\smallfrown="3\msa@61
\mathchardef\Subset="3\msa@62
\mathchardef\Supset="3\msa@63
\mathchardef\Cup="2\msa@64

\mathchardef\Cap="2\msa@65

\mathchardef\curlywedge="2\msa@66
\mathchardef\curlyvee="2\msa@67
\mathchardef\leftthreetimes="2\msa@68
\mathchardef\rightthreetimes="2\msa@69
\mathchardef\subseteqq="3\msa@6A
\mathchardef\supseteqq="3\msa@6B
\mathchardef\bumpeq="3\msa@6C
\mathchardef\Bumpeq="3\msa@6D
\mathchardef\lll="3\msa@6E

\mathchardef\ggg="3\msa@6F

\mathchardef\circledS="0\msa@73
\mathchardef\pitchfork="3\msa@74
\mathchardef\dotplus="2\msa@75
\mathchardef\backsim="3\msa@76
\mathchardef\backsimeq="3\msa@77
\mathchardef\complement="0\msa@7B
\mathchardef\intercal="2\msa@7C
\mathchardef\circledcirc="2\msa@7D
\mathchardef\circledast="2\msa@7E
\mathchardef\circleddash="2\msa@7F
\def\ulcorner{\delimiter"4\msa@70\msa@70 }
\def\urcorner{\delimiter"5\msa@71\msa@71 }
\def\llcorner{\delimiter"4\msa@78\msa@78 }
\def\lrcorner{\delimiter"5\msa@79\msa@79 }
\def\yen{\mathhexbox\msa@55 }
\def\checkmark{\mathhexbox\msa@58 }
\def\circledR{\mathhexbox\msa@72 }
\def\maltese{\mathhexbox\msa@7A }
\mathchardef\lvertneqq="3\msb@00
\mathchardef\gvertneqq="3\msb@01
\mathchardef\nleq="3\msb@02
\mathchardef\ngeq="3\msb@03
\mathchardef\nless="3\msb@04
\mathchardef\ngtr="3\msb@05
\mathchardef\nprec="3\msb@06
\mathchardef\nsucc="3\msb@07
\mathchardef\lneqq="3\msb@08
\mathchardef\gneqq="3\msb@09
\mathchardef\nleqslant="3\msb@0A
\mathchardef\ngeqslant="3\msb@0B
\mathchardef\lneq="3\msb@0C
\mathchardef\gneq="3\msb@0D
\mathchardef\npreceq="3\msb@0E
\mathchardef\nsucceq="3\msb@0F
\mathchardef\precnsim="3\msb@10
\mathchardef\succnsim="3\msb@11
\mathchardef\lnsim="3\msb@12
\mathchardef\gnsim="3\msb@13
\mathchardef\nleqq="3\msb@14
\mathchardef\ngeqq="3\msb@15
\mathchardef\precneqq="3\msb@16
\mathchardef\succneqq="3\msb@17
\mathchardef\precnapprox="3\msb@18
\mathchardef\succnapprox="3\msb@19
\mathchardef\lnapprox="3\msb@1A
\mathchardef\gnapprox="3\msb@1B
\mathchardef\nsim="3\msb@1C
\mathchardef\napprox="3\msb@1D
\mathchardef\nsubseteqq="3\msb@22
\mathchardef\nsupseteqq="3\msb@23
\mathchardef\subsetneqq="3\msb@24
\mathchardef\supsetneqq="3\msb@25
\mathchardef\subsetneq="3\msb@28
\mathchardef\supsetneq="3\msb@29
\mathchardef\nsubseteq="3\msb@2A
\mathchardef\nsupseteq="3\msb@2B
\mathchardef\nparallel="3\msb@2C
\mathchardef\nmid="3\msb@2D
\mathchardef\nshortmid="3\msb@2E
\mathchardef\nshortparallel="3\msb@2F
\mathchardef\nvdash="3\msb@30
\mathchardef\nVdash="3\msb@31
\mathchardef\nvDash="3\msb@32
\mathchardef\nVDash="3\msb@33
\mathchardef\ntrianglerighteq="3\msb@34
\mathchardef\ntrianglelefteq="3\msb@35
\mathchardef\ntriangleleft="3\msb@36
\mathchardef\ntriangleright="3\msb@37
\mathchardef\nleftarrow="3\msb@38
\mathchardef\nrightarrow="3\msb@39
\mathchardef\nLeftarrow="3\msb@3A
\mathchardef\nRightarrow="3\msb@3B
\mathchardef\nLeftrightarrow="3\msb@3C
\mathchardef\nleftrightarrow="3\msb@3D
\mathchardef\divideontimes="2\msb@3E
\mathchardef\varnothing="0\msb@3F
\mathchardef\nexists="0\msb@40
\mathchardef\mho="0\msb@66
\mathchardef\thorn="0\msb@67
\mathchardef\beth="0\msb@69
\mathchardef\gimel="0\msb@6A
\mathchardef\daleth="0\msb@6B
\mathchardef\lessdot="3\msb@6C
\mathchardef\gtrdot="3\msb@6D
\mathchardef\ltimes="2\msb@6E
\mathchardef\rtimes="2\msb@6F
\mathchardef\shortmid="3\msb@70
\mathchardef\shortparallel="3\msb@71
\mathchardef\smallsetminus="2\msb@72
\mathchardef\thicksim="3\msb@73
\mathchardef\thickapprox="3\msb@74
\mathchardef\approxeq="3\msb@75
\mathchardef\succapprox="3\msb@76
\mathchardef\precapprox="3\msb@77
\mathchardef\curvearrowleft="3\msb@78
\mathchardef\curvearrowright="3\msb@79
\mathchardef\digamma="0\msb@7A
\mathchardef\varkappa="0\msb@7B
\mathchardef\hslash="0\msb@7D
\mathchardef\hbar="0\msb@7E
\mathchardef\backepsilon="3\msb@7F
\def\Bbb{\ifmmode\let\next\Bbb@\else
 \def\next{\errmessage{Use \string\Bbb\space only in math mode}}\fi\next}
\def\Bbb@#1{{\Bbb@@{#1}}}
\def\Bbb@@#1{\fam\msbfam#1}

\catcode`\@=\active

 \def\CS{\hbox{{$\cal S$}}}

\def\CR{\hbox{{$\cal R$}}} 
 
\def\CM{\hbox{{$\cal M$}}}
\def\CN{\hbox{{$\cal N$}}}

\def\CO{\hbox{{$\cal O$}}}

\def\R{{\Bbb R}}
\def\C{{\Bbb C}}

\def\cu{{\sl u}}

\def\del{{\partial}}
\def\lform{\hbox{$\sqcup$}\llap{\hbox{$\sqcap$}}}

\def\eps{{\epsilon}}

\def\trace{{\rm Tr\, }}

\def\<{\langle}
\def\>{\rangle}
\def\rcross{{\triangleright\!\!\!<}}
\def\lcross{{>\!\!\!\triangleleft}}

\def\cobicross{{\triangleright\!\!\!\blacktriangleleft}}
\def\bicross{{\blacktriangleright\!\!\!\triangleleft}}
\def\dcross{{\bowtie}}

\def\tens{\mathop{\otimes}}
\def\la{{\triangleright}}\def\ra{{\triangleleft}}
\def\isom{{\cong}}

\def\image{{\rm image}\, }

\def\id{{\rm id}}

\def\o{{}_{(1)}}\def\t{{}_{(2)}}

\def\uo{{{}^{(1)}}}\def\ut{{{}^{(2)}}}

\def\extd{{\rm d}}

\def\proof{\goodbreak\noindent{\bf Proof\quad}}
\def\endproof{{\ $\lform$}\bigskip }

\def\text#1{{\rm #1}}
\def\note#1{}

\def\nquad{{\!\!\!\!\!\!}}

\def\equad{\nquad}
\def\eqn#1#2{\begin{equation}#2\label{#1}\end{equation}}

\def\align#1{\begin{eqnarray*}#1\end{eqnarray*}}

\def\cmath#1{\[\begin{array}{c} #1 \end{array}\]}

\def\ceqn#1#2{\begin{equation}\label{#1}\begin{array}{c}#2\end{array}
\end{equation}}

\documentstyle[11pt]{article}
\newtheorem{lemma}{Lemma}[section]
\newtheorem{propos}[lemma]{Proposition}
\newtheorem{example}[lemma]{Example}
\newtheorem{theorem}[lemma]{Theorem}

\newtheorem{corol}[lemma]{Corollary}

\begin{document}
\baselineskip 22pt

\newpage
\baselineskip 22pt

{\ }\hskip 4.7in   Damtp/96-97
\vspace{.2in}

\begin{center} {\Large QUASITRIANGULAR AND DIFFERENTIAL STRUCTURES ON\\
BICROSSPRODUCT HOPF ALGEBRAS}
\\ \baselineskip 13pt{\ }
{\ }\\
Edwin Beggs \\
{\ }\\
Department of Mathematics\\
University of Wales, Swansea\\
Singleton Park, Swansea SA2 8PP, UK\\
{\ }\\
+\\
{\ }\\
Shahn Majid\footnote{Royal Society University
Research Fellow and Fellow of Pembroke College, Cambridge}\\
{\ }\\
Department of Mathematics, Harvard University\\
Science Center, Cambridge MA 02138, USA\footnote{During 1995 \& 1996}\\
+\\
Department of Applied Mathematics and Theoretical Physics\\
University of Cambridge, Cambridge CB3 9EW, UK
\end{center}

\begin{center}
December 1996
\end{center}
\vspace{10pt}
\begin{quote}\baselineskip 13pt
\noindent{\bf ABSTRACT} Let $X=GM$ be a finite group factorisation. It is shown
that the quantum double
$D(H)$ of the associated bicrossproduct Hopf algebra $H=kM\cobicross k(G)$ is
itself a bicrossproduct $kX\cobicross k(Y)$ associated to a group $YX$, where
$Y=G\times M^{\text{op}}$. This
provides a  class of bicrossproduct Hopf algebras which are
quasitriangular. We also construct a subgroup $Y^\theta X^\theta$
associated to every   order-reversing automorphism $\theta$ of $X$. The
corresponding Hopf algebra  $kX^\theta\cobicross k(Y^\theta)$ has the same
coalgebra as $H$. Using related results,
we classify the first order bicovariant differential calculi on $H$ in terms of
orbits in a certain quotient space of $X$.
\end{quote}
\baselineskip 22pt

\section{INTRODUCTION}

The quantum double\cite{Dri}   of the bicrossproduct Hopf algebra
$H=kM\cobicross k(G)$ associated to a finite group factorisation $X=GM$ has
been studied recently in \cite{BGM:fin}. Here we continue this study with
further results in the same topic, including a concrete application to the
classification of the bicovariant differential calculi on a bicrossproduct.

The bicrossproduct Hopf algebras have been introduced in \cite{Tak:mat} and
\cite{Ma:phy}, and extensively studied since then. Factorisations of groups
abound in mathematics, so these Hopf algebras, which are non-commutative and
non-cocommutative, are quite common. In the context of \cite{Ma:phy} they are
viewed as systems which combine quantum mechanical ideas with geometry in a
unified way. To develop this idea a natural step is to compute and study the
algebra of differential forms (the so-called differential calculus) over them,
which we do here. Another open problem is the close connection which one may
expect between these Hopf algebras and the method of inverse scattering in
soliton theory. We make a connection of this type also in the paper. All our
results, however, will be in a strictly algebraic setting with finite groups
and finite-dimensional Hopf algebras.

We recall two principal results from \cite{BGM:fin}. One result was that $D(H)$
is isomorphic to a twist (i.e. up to a categorical equivalence of
representations) to the much simpler double $D(X)$ of the group algebra $kX$.
We provide in Section~2 a variant of this result, namely that $D(H)$ is also
isomorphic to a bicrossproduct $k X\cobicross k(Y)$ where $Y$ is the `dressing
group' $G\times M^{\rm op}$ (the terminology comes from soliton
theory\cite{Sem:dre}). The
bicrossproduct is associated to a certain double cross product group $Y\dcross
X$ which factorises into $X,Y$. Our new result answers affirmatively the
question
`can a bicrossproduct Hopf algebra be quasitriangular?'. Among non-commutative
and non-cocommutative Hopf algebras, the quasitriangular\cite{Dri} ones have a
special place with special properties. We compute the quasitriangular structure
in the present case, and some of its consequences.

The second principal result in \cite{BGM:fin} was that
associated to  every  order reversing isomorphism $\theta$ of
$X$ is a Hopf algebra isomorphism $\Theta$ of the quantum double. In Section~3
we provide more results about $\Theta$. We then construct two
new groups $Y^\theta,X^\theta$ forming a
subgroup $Y^\theta \dcross X^\theta$ of $Y\dcross X$.
We study the associated Hopf algebra $kX^\theta\cobicross k(Y^\theta)$. It
has an isomorphic coalgebra to that of our original bicrossproduct $H$
associated to  the
factorisation $GM$. However we give a
finite group example where these are not isomorphic
as algebras.
Another example involving upper and lower
triangular matrices is given. In this case we use the
order reversing isomorphism given by the operation
$g\mapsto (g^*)^{-1}$.

In Section~4 we study a certain Hilbert space representation of $H$ and find
that our various constructions over $\C$ respect the $*$-structures. This is
also
 one of the motivations behind the main isomorphism in Section~2.

In Section~5 we turn to a specific application of the quantum double, namely to
the construction of first order bicovariant differential calculi. By definition
a first order differential calculus over an algebra $A$ means an $A$-bimodule
$\Omega^1$ over $A$ and a map $\extd:A\to\Omega^1$ obeying the Leibniz rule
(which
makes sense using the bimodule structure). When $A$ is a Hopf algebra the
calculus is said to be left-, right- or bi-covariant when $\extd$ interwines
the
left regular coaction of $A$, the right regular coaction, or both, with a
coaction given on $\Omega^1$. It is known that first order bicovariant calculi
are related to the representation theory of the quantum double\cite{Wor:dif},
allowing us to apply our previous results\cite{BGM:fin} when $A=k(M)\cobicross
kG$.
We find that the irreducible
bicovariant calculi correspond to the choice of an orbit in
a certain quotient space of the group $X$ along with an irreducible
subrepresentation of the isotropy subgroup associated to the orbit. We actually
classify the bicovariant quantum tangent spaces using the techniques developed
in \cite{Ma:cla}, and obtain the corresponding 1-forms later by dualisation.
The result is a constructive method which provides the entire moduli space of
bicovariant calculi on a bicrossproduct, as we demonstrate on some nontrivial
examples.

We note that in the
physics literature an important bicrossproduct Hopf algebra is the
$\kappa$-deformed
Poincar\'e algebra\cite{MaRue:bic}, and for this Hopf algebra some examples of
bicovariant
calculi have been obtained by other means\cite{KMS:bic}. Here the group $G$ is
the
Lorentz group and $M$ is $\R^3\lcross\R$. Another important (and very similar)
example is with $G=SO(3)$ and $M=\R^2\lcross\R$ in \cite{Ma:phy}\cite{Ma:book},
where the bicrossproduct is the algebra of observables (or quantum phase space)
of a particular quantum mechanical system.
It seems likely that a more geometrical version of the present results should
include such examples as well. Armed with a choice of differential calculus,
one may proceed to gauge theory (i.e. to bundles and connections) using the
formalism of \cite{BrzMa:gau}.

\section*{Preliminaries}

We use the notation and conventions of \cite{BGM:fin}, which are also the
notations and
conventions in the text \cite{Ma:book}. Briefly, let $X=GM$ be a group which
factorises into two subgroups. Then each group
acts on the other through left and right actions $\la,\ra$ defined by $su=(s\la
u)(s\ra u)$ forall $s\in M$ and $u\in G$. Conversely, given two actions
$\la,\ra$ obeying certain matching conditions

\ceqn{mpair}
{s\ra e=s,\quad (s\ra u)\ra v=s\ra(uv);\qquad e\ra u=e,\quad
(st)\ra
u=\left(s\ra(t\la u)\right)(t\ra u) \\
e\la u=u,\quad s\la(t\la u)=(st)\la u;\quad s\la e=e,\quad s\la(uv)=(s\la
u)\left((s\ra u)\la v\right).}

we can build a double cross product group on $G\times M$ with
\eqn{dgroup}{(u,s)(v,t)=(u(s\la v),(s\ra v)t),\quad e=(e,e),\quad
(u,s)^{-1}=(s^{-1}\la
u^{-1},s^{-1}\ra u^{-1}).}
The associated bicrossproduct Hopf algebra $H=k M\cobicross k(G)$ has the smash
product algebra structure by the induced action of $M$ and the smash coproduct
coalgebra structure by the induced coaction of $G$. Explicitly,
\ceqn{bic}{(s\tens \delta_u)(t\tens \delta_v)=\delta_{u,t\la v}(st\tens
\delta_{v}),\quad \Delta (s\tens \delta_u)=\sum_{xy=u}s\tens\delta_x\tens
s\ra x\tens \delta_y\\
1=\sum_u e\tens \delta_u,\quad \eps (s\tens \delta_u)=\delta_{u,e},\quad
S(s\tens \delta_u)=(s\ra u)^{-1}\tens\delta_{(s\la u)^{-1}}.}
We work over a general ground field $k$. There is also a natural
$*$-algebra structure $(s\tens \delta_u)^*=s^{-1}\tens\delta_{s\la u}$ when the
ground field has an involution. This happens over $\C$, but can also be
supposed for any field with $\bar\lambda=\lambda$ for all $\lambda\in k$. The
dual $H^*$ has a similar structure $k(M)\cobicross kG$ on the dual basis,
\ceqn{cobic}{(\delta_s\tens u)(\delta_t\tens v)=\delta_{s\ra u,t}(\delta_s\tens
uv),
\quad
\Delta(\delta_s\tens u)=\sum_{ab=s}\delta_a\tens b\la u\tens
\delta_b\tens u\\
1=\sum_s\delta_s\tens e,\quad \eps(\delta_s\tens u)=\delta_{s,e},\quad
S(\delta_s\tens u)=\delta_{(s\ra u)^{-1}}\tens(s\la u)^{-1},}
and $(\delta_s\tens u)^*=\delta_{s\ra u}\tens u^{-1}$ when the ground field has
an involution. The quantum
double\cite{Dri} is a general construction $D(H)=H^{*\rm op}\dcross H$ built on
$H^*\tens H$ with a double cross product algebra structure and tensor product
coalgebra structure. In our case the cross relations between $H,H^{*\rm op}$
are\cite{BGM:fin}
\ceqn{double}{(1\tens t\tens\delta_v)(\delta_s\tens u\tens 1)=\delta_{t's(t\ra
vu^{-1})^{-1}}\tens (t\ra vu^{-1})\la u\tens t'\tens \delta_{(s\la u)vu^{-1}},}
where $t'=t\ra (s\la u)^{-1}$.

\section{More about $D(H)$}

Here we extend results about the quantum double associated to a bicrossproduct
in \cite{BGM:fin}.
For our first observation, it is known that to every factorisation $X=GM$ there
is a `double factorisation' $YX$ where $Y$ is also $G\times M$ as a set and the
action of $X$  is the adjoint action viewed as an action on $Y$\cite{Ma:book}.
Here we give a similar but different `double factorisation' more suitable for
our needs.

\begin{propos} Let $Y=G\times M^{\text{op}}$ with group law
$(us).(vt)=uvts$. Then there is a double cross product group $Y\dcross X$
(factorising into $Y,X$) defined by actions
\[ us\tilde\ra vt=((s\ra v)ts^{-1}\la u^{-1})^{-1}(s\ra v)\]
\[ us\tilde\la vt=us(vt)(us)^{-1}=u(s\la v)((s\ra v)ts^{-1}\la u^{-1})((s\ra
v)ts^{-1}\ra u^{-1}).\]
The second line is the adjoint action on $X$ which we view as an action on the
set $Y$.
\end{propos}
\proof  We show that these actions are matched in the required sense (see
\cite{Ma:book}). First note that
the following results are immediate from the definitions:
\[
us\tilde\ra e=us\ ,\ e\tilde\ra vt=e\ ,\
us\tilde\la e=e\ ,\ e\tilde\la vt=vt\ .
\]
For the other results we must work rather harder. To derive the equation
\[
us\tilde\la \big( (vt).(wr) \big)\ =\
\big(us\tilde\la vt\big).
\big((us\tilde\ra vt)\tilde\la wr\big)\ ,
\]
begin by the formula given for $\tilde \la$ as
$
us\tilde\la vt=u(s\la v)yp
$
where
$
y = (s\ra v)ts^{-1}\la u^{-1}$ and
$p = (s\ra
v)ts^{-1}\ra u^{-1}$.
Then starting from $us\tilde\ra vt=y^{-1}(s\ra v)$
 we calculate
\[
(us\tilde\ra vt)\tilde\la wr\ =\ y^{-1}\,
\big((s\ra v)\la w\big)\,
\big((s\ra vw)r(s\ra v)^{-1}\la y\big)\,
\big((s\ra vw)r(s\ra v)^{-1}\ra y\big)\ ,
\]
and on applying the rules for multiplication in $Y$
we find the required formula,
\begin{eqnarray*}
\big(us\tilde\la vt\big).
\big((us\tilde\ra vt)\tilde\la wr\big) &=&
u\, (s\la v)\, \big((s\ra v)\la w\big)\,
\big((s\ra vw)r(s\ra v)^{-1}\la y\big)\,
\big((s\ra vw)r(s\ra v)^{-1}\ra y\big)\,  p\\
 &=&
u\, (s\la vw)\,
\big((s\ra vw)rts^{-1}\la u^{-1}\big)\,
\big((s\ra vw)rts^{-1}\ra u^{-1}\big)\\
 &=&
us\tilde\la \big( (vt).(wr) \big)
\ .
\end{eqnarray*}
Now we must prove that
\[
\big( (wr).(us)\big)\tilde\ra vt\ =\
\big(wr\tilde\ra (us\tilde\la vt)\big).
\big(us\tilde\ra vt\big)\ .
\]
Begin by calculating
\[wr\tilde\ra (us\tilde\la vt)\ =\
wr\tilde\ra (u(s\la v)yp)\ =\
\big((r\ra u(s\la v)y) pr^{-1}\la w^{-1}\big)^{-1}
\, \big(r\ra u(s\la v)y\big)\ ,
\]
and then use the definition of multiplication in $X$ to find
\begin{eqnarray*}
\big(wr\tilde\ra (us\tilde\la vt)\big).
\big(us\tilde\ra vt\big)& =&
\big((r\ra u(s\la v)y) pr^{-1}\la w^{-1}\big)^{-1}
\, \big(r\ra u(s\la v)y\big)\, y^{-1}\, (s\ra v)  \\
& =&
\big((r\ra u(s\la v)y) pr^{-1}\la w^{-1}\big)^{-1}
\, \big(\big(r\ra u(s\la v)y\big)\la y^{-1}\big)\,
 \big(r\ra u(s\la v)\big)\, (s\ra v)
 \\
& =&
\big((r\ra u(s\la v)y) pr^{-1}\la w^{-1}\big)^{-1}
\, \big(\big(r\ra u(s\la v)\big)\la y\big)^{-1}\,
 \big((r\ra u)s\ra v\big)
 \\
& =&
\big[ \big(\big(r\ra u(s\la v)\big)\la y\big) \,
\big((r\ra u(s\la v)y) pr^{-1}\la w^{-1}\big)
\big]^{-1}\,
 \big((r\ra u)s\ra v\big)
 \\
& =&
\big[ \big(r\ra u(s\la v)\big)\la y
( pr^{-1}\la w^{-1})
\big]^{-1}\,
 \big((r\ra u)s\ra v\big)
\ .
\end{eqnarray*}
Now consider a part of this equation
\begin{eqnarray*}
y( pr^{-1}\la w^{-1})& =&
\big((s\ra v)ts^{-1}\la u^{-1}\big)\,
\big(((s\ra
v)ts^{-1}\ra u^{-1})\la(r^{-1}\la w^{-1})\big)  \\
& =&
(s\ra v)ts^{-1}\la u^{-1}(r^{-1}\la w^{-1})\ ,
\end{eqnarray*}
and putting this into the previous equation we get
\begin{eqnarray*}
\big(wr\tilde\ra (us\tilde\la vt)\big).
\big(us\tilde\ra vt\big)
& =&
\big[ \big(r\ra u(s\la v)\big)(s\ra v)ts^{-1}\la u^{-1}(r^{-1}\la w^{-1})
\big]^{-1}\,
 \big((r\ra u)s\ra v\big)   \\
& =&
\big[ \big((r\ra u)s\ra v\big)ts^{-1}\la u^{-1}(r^{-1}\la w^{-1})
\big]^{-1}\,
 \big((r\ra u)s\ra v\big)
\ .
\end{eqnarray*}
On comparing this with
\[
\big( (wr).(us)\big)\tilde\ra vt\ =\
\big( w(r\la u)(r\la u)s)\big)\tilde\ra vt\ =\
\big[
\big((r\ra u)s\ra v\big) t s^{-1}(r\ra u)^{-1} \la
(r\la u)^{-1}w^{-1}
\big]^{-1}
\, \big((r\ra u)s\ra v\big)
\]
we get identical results since
\[
(r\ra u)^{-1} \la
(r\la u)^{-1}w^{-1}\ =\
\big((r\ra u)^{-1} \la
(r\la u)^{-1}\big)\, \big(\big((r\ra u)^{-1} \ra
(r\la u)^{-1}\big)\la w^{-1}\big)\ =\
u^{-1}(r^{-1}\la w^{-1})
\]
 \endproof

\begin{theorem}\label{Dbic} $D(H)\isom k X \cobicross k(Y)$ as Hopf algebras,
by
\[\psi:D(H)\to k X\cobicross k(Y),\quad \psi(\delta_s\tens u\tens
t\tens\delta_v)= (s\la u)^{-1}t\tens\delta_{v(t\ra
v)^{-1}s^{-1}t}\ .\]
Over a field with involution, the map preserves the star operation.
\end{theorem}
\proof  The structure of $D(H)$ in the basis used is in \cite{BGM:fin}. We
check that the linear map $\psi$ is an algebra isomorphism
 to the smash product induced by $\tilde\la$.
Start with
 $\alpha=\delta_x\tens q\tens t\tens\delta_v$, and
$\beta=\delta_s\tens u\tens r\tens\delta_w$ in $D(H)$, and
multiply them together to get
\[
\alpha \beta\ =\
\delta_{s'\ra u',x}\ \delta_{v',r\la w}\ \big(
\delta_{s'} \tens u'q \tens t'r \tens \delta_{w} \big)\ ,
\]
where
\[
t'=t\ra(s\la u)^{-1}\ ,\ v'=(s\la u) vu^{-1}\ ,\ s'=t's(t\ra vu^{-1})^{-1}\ ,\
u'=(t\ra vu^{-1})\la u\ .
\]
Now we can calculate
\begin{eqnarray*}
\psi(\alpha\beta)& =&
\delta_{s'\ra u',x}\ \delta_{v',r\la w}\ \big(
(s'\la u'q)^{-1}t'r\tens \delta_{w(t'r\ra w)^{-1}s'^{-1}t'r}\big) \\
& =&
\delta_{s'\ra u',x}\ \delta_{v',r\la w}\ \big(
((s'\ra u')\la q)^{-1}  (s'\la u')^{-1}
t'r\tens \delta_{w
(r\ra w)^{-1}
(t'\ra(r\la w))^{-1}
s'^{-1}t'r}\big) \\
& =&
\delta_{t(s\ra u)(t\ra v)^{-1},x}\ \delta_{
(s\la u) vu^{-1},r\la w}\ \big(
(x\la q)^{-1}  (t's\la u)^{-1}
t'r\tens \delta_{w
(r\ra w)^{-1}
(t'\ra v')^{-1}
s'^{-1}t'r}\big)
 \\
& =&
\delta_{t(s\ra u)(t\ra v)^{-1},x}\ \delta_{
(s\la u) vu^{-1},r\la w}\ \big(
(x\la q)^{-1}
(t'\ra(s\la u))
(s\la u)^{-1}
r\tens \delta_{w
(r\ra w)^{-1}
s^{-1}r}\big)
\\
& =&
\delta_{t(s\ra u)(t\ra v)^{-1},x}\ \delta_{
(s\la u) vu^{-1},r\la w}\ \big(
(x\la q)^{-1}
t
(s\la u)^{-1}
r\tens \delta_{w
(r\ra w)^{-1}
s^{-1}r}\big)
\ .
\end{eqnarray*}
Here we have used
\begin{eqnarray*}
s'\ra u'& = &\big(t\ra (s\la u)^{-1}\big)s(t\ra vu^{-1})^{-1}    \ra
\big((t\ra vu^{-1})\la u\big)  \\
& = & \big(  \big(t\ra (s\la u)^{-1}\big)s \ra u \big)
 \big( (t\ra vu^{-1})^{-1} \ra  \big( (t\ra vu^{-1}) \la u \big) \big) \\ & = &
 \big(t\ra (s\la u)^{-1}\ra (s\la u) \big) (s\ra u)(t\ra vu^{-1}\ra u)^{-1} \\
& = & t(s\ra u)(t\ra v)^{-1} \ .
\end{eqnarray*}
Conversely we can calculate the product in $k X \cobicross k(Y)$ as
\begin{eqnarray*}
\psi(\alpha)\psi(\beta)& =&
\big(
 (x\la q)^{-1} t \tens \delta_{v  (t\ra v)^{-1}x^{-1} t  } \big)
\big(
 (s\la u)^{-1} r\tens \delta_{w  (r\ra w)^{-1} s^{-1} r  } \big)  \\
& =& \delta_{v  (t\ra v)^{-1}x^{-1} t\, ,\,
(s\la u)^{-1} r \tilde\la w  (r\ra w)^{-1} s^{-1} r  }\
 (x\la q)^{-1} t
 (s\la u)^{-1} r\tens \delta_{w  (r\ra w)^{-1} s^{-1} r  }
 \\
& =& \delta_{v  (t\ra v)^{-1}x^{-1} t\, ,\,
(s\la u)^{-1} (r \la w)  u(s\ra u)^{-1}  }\
 (x\la q)^{-1} t
 (s\la u)^{-1} r\tens \delta_{w  (r\ra w)^{-1} s^{-1} r  }
 \\
& =& \delta_{v  \, ,\,
(s\la u)^{-1} (r \la w)  u  }\
\delta_{  (t\ra v)^{-1}x^{-1} t\, ,\,
(s\ra u)^{-1}  }\
 (x\la q)^{-1} t
 (s\la u)^{-1} r\tens \delta_{w  (r\ra w)^{-1} s^{-1} r  }\ .
\end{eqnarray*}
It is now apparent that the expressions
for $\psi(\alpha\beta)$ and $\psi(\alpha)\psi(\beta)$
are the same.

To compare the coproducts, we use the coproduct of $D(H)$, which is the tensor
product one
\[
\Delta_{D(H)}\big(\delta_s \tens u\tens t\tens\delta_v\big) \
=\
\sum_{ab=s\ ,\ xy=v}
(\delta_a \tens b\la u\tens t\tens\delta_x)
\tens
(\delta_b \tens u\tens t\ra x\tens\delta_y)
\]
Applying $\psi\tens\psi$ to this, we
find the following expression for $
(\psi\tens\psi)
\Delta_{D(H)}\big(\delta_s \tens u\tens t\tens\delta_v\big)$ :
\[
\sum_{ab=s\ ,\ xw=v}
\big((s\la u)^{-1}t \tens \delta_{x(t\ra x)^{-1}a^{-1}t}\big)
\tens
\big((b\la u)^{-1}(t\ra x) \tens \delta_{w(t\ra v)^{-1}
b^{-1}(t\ra x)}\big)
\ ,
\]
Alternatively we can calculate
\begin{eqnarray*}
\Delta_{kX\cobicross k(Y)}\psi
(\delta_s \tens u\tens t\tens\delta_v)& =&
\Delta_{kX\cobicross k(Y)}\big(
(s\la u)^{-1}t \tens \delta_{v(t\ra v)^{-1}s^{-1}t}
\big)   \\ & =&
\sum_{yz=v(t\ra v)^{-1}s^{-1}t}
(s\la u)^{-1}t \tens \delta_{y} \tens
(s\la u)^{-1}t\tilde\ra y \tens \delta_{z} \ ,
\end{eqnarray*}
where $y,z\in Y$,
which is the smash coproduct for the coaction induced by the back-reaction
$\tilde\ra$. We begin with the calculation
\[
\big(x(t\ra x)^{-1}a^{-1}t\big).\big(w(t\ra v)^{-1}
b^{-1}(t\ra x)\big)\ =\
xw(t\ra v)^{-1} b^{-1}a^{-1}t\ ,
\]
which shows that if we replace $y$ by $x(t\ra x)^{-1}a^{-1}t$ and $z$ by
$w(t\ra v)^{-1}b^{-1}(t\ra x)$, that the conditions of the summations
are the same. It now remains to calculate
\[
(s\la u)^{-1}t\tilde\ra y \ = \ (s\la u)^{-1}t\tilde\ra
x(t\ra x)^{-1}a^{-1}t\ =\
(a^{-1}s\la u)^{-1}(t\ra x)\ =\ (b\la u)^{-1}(t\ra x)\ ,
\]
as required.

If we apply $\psi$ to the unit of $D(H)$, we get
\[
\psi \sum_{s,v} \delta_s\tens e\tens e\tens\delta_v\ =\
\sum_{s,v}  e\tens \delta_{vs^{-1}}\ ,
\]
which is the unit in $kX \cobicross k(Y)$.

Next we consider the counit,
\[
\eps_{k X \cobicross k(Y)}\big(\psi(
\delta_s\tens u\tens t\tens\delta_v)\big)\ =\
\eps_{k X \cobicross k(Y)}\big(
(s\la u)^{-1}t\tens\delta_{v(t\ra
v)^{-1}s^{-1}t}\big)\ =\ \delta_{v(t\ra
v)^{-1}s^{-1}t,e}
\]
The last $\delta$-function splits into $\delta_{v,e}
\delta_{(t\ra v)^{-1}s^{-1}t,e}$, which is equal to
$\delta_{v,e}\delta_{s,e}$, which is in turn equal to
$\eps_{D(H)}\big(\delta_s\tens u\tens t\tens\delta_v\big)$.

Finally we consider the antipode,
\begin{eqnarray*}
S_{k X \cobicross k(Y)}\big(\psi(
\delta_s\tens u\tens t\tens\delta_v)\big)& =&
S_{k X \cobicross k(Y)}\big(
(s\la u)^{-1}t\tens\delta_{v(t\ra v)^{-1}s^{-1}t}\big)   \\
& =&
\big( (s\la u)^{-1}t \tilde \ra v(t\ra v)^{-1}s^{-1}t \big)^{-1} \tens
\delta_{\big( (s\la u)^{-1}t \tilde \la v(t\ra v)^{-1}s^{-1}t \big)^{-1}}
 \\
& =&
\big( u^{-1}(t\ra v) \big)^{-1} \tens
\delta_{\big( (s\la u)^{-1}(t\la v) s^{-1}(s\la u) \big)^{-1}} \\
& =&
(t\ra v)^{-1}u  \tens
\delta_{\big( (s\la u)^{-1}(t\la v)u(s\ra u)^{-1} \big)^{-1}}\\
& =&
(t\ra v)^{-1}u  \tens
\delta_{ u^{-1}(t\la v)^{-1}(s\la u)(s\ra u) } \ ,
\end{eqnarray*}
where we remember in the last line to take the inverse for the $Y$
group operation,
and compare this with
\begin{eqnarray*}
\psi S_{D(H)}\big( \delta_s\tens u\tens t\tens\delta_v \big) & = &
\psi\big((1\tens S(t\tens\delta_v))
(S^{-1}(\delta_s\tens u)\tens 1)\big)  \\
& = &
\psi\big((1\tens (t\ra v)^{-1}\tens\delta_{(t\la v)^{-1}})
(\delta_{(s\ra u)^{-1}}\tens (s\la u)^{-1}\tens 1)\big)
 \\
& = &
\psi\big(\delta_{t'\bar s(\bar t\ra
\bar v\bar u^{-1})^{-1}}
\tens (\bar t\ra \bar v\bar u^{-1})\la \bar u\tens t'\tens
\delta_{(\bar s\la \bar u)\bar v\bar u^{-1}}
\big)
\end{eqnarray*}
where $\bar t=(t\ra v)^{-1}$,
$\bar v=(t\la v)^{-1}$,
$\bar s=(s\ra u)^{-1}$,
$\bar u=(s\la u)^{-1}$ and $t'=\bar t\ra (\bar s\la \bar u)^{-1}$.
Applying the definition of $\psi$, and using the fact that
$\bar s\la \bar u=u^{-1}$, we find
\begin{eqnarray*}
\psi S_{D(H)}\big( \delta_s\tens u\tens t\tens\delta_v \big) & =&
(t'\bar s\la \bar u)^{-1}t'\, \tens\, \delta_{
(\bar s\la \bar u)\bar v\bar u^{-1}
(t'\ra(\bar s\la \bar u)\bar v\bar u^{-1})^{-1}
(\bar t\ra \bar v\bar u^{-1}) \bar s^{-1}t'^{-1}t'}   \\
& =&
(t'\la  u^{-1})^{-1}t'\, \tens\, \delta_{
u^{-1}\bar v\bar u^{-1}
 \bar s^{-1} }   \\
& =&
(t'\ra  u^{-1}) u\, \tens\, \delta_{ u^{-1}
\bar v (s\la u)(s\ra u)}    \\
& =&
(t\ra  v)^{-1} u\, \tens\, \delta_{ u^{-1}
(t\la v)^{-1}su }\ ,
\end{eqnarray*}
again as required. This concludes the proof of the Hopf algebra
isomorphism. Now we show that the star operation
is preserved.
\begin{eqnarray*}
\psi * \big( \delta_s\tens u\tens t\tens\delta_v \big) & = &
\psi\big( (I\tens t^{-1}\tens \delta_{t\la v})
\ (\delta_{s\ra u}\tens u^{-1}\tens I) \big)  \\
& = &
\psi\big(
\delta_{t'(s\ra u)(t^{-1}\ra(t\la v)u)^{-1}}  \tens (t^{-1}\ra(t\la v)u)\la
u^{-1}
\tens t' \tens \delta_{((s\ra u)\la u^{-1})(t\la v)u} \big)  \ ,
\end{eqnarray*}
where $t'=t^{-1}\ra (s\la u)$. Applying the definition of $\psi$, we get
\begin{eqnarray*}
\psi * \big( \delta_s\tens u\tens t\tens\delta_v \big) & = &
(t'(s\ra u)\la u^{-1})^{-1}t' \ \tens\
\delta_{(s\la u)^{-1}(t\la v) u \big(t^{-1}\ra(t\la v)u\big)^{-1}
 \big(t^{-1}\ra(t\la v)u\big)(s\ra u)^{-1}t'^{-1}t' } \\
& = & ((t^{-1}s\ra u)\la u^{-1})^{-1}t' \ \tens\ \delta_{
(s\la u)^{-1}(t\la v) u (s\ra u)^{-1}}   \\
& = & t^{-1}(s\la u) \ \tens\ \delta_{
(s\la u)^{-1}(t\la v) u (s\ra u)^{-1}}\ ,   \\
\phantom{y} * \psi  \big( \delta_s\tens u\tens t\tens\delta_v \big) & = &
*\big((s\la u)^{-1}t\tens\delta_{v(t\ra
v)^{-1}s^{-1}t}\big)  \\
& = & t^{-1}(s\la u) \tens\delta_{(s\la u)^{-1}tv(t\ra
v)^{-1}s^{-1}tt^{-1}(s\la u)}  \\
& = & t^{-1}(s\la u) \tens\delta_{(s\la u)^{-1}(t\la v) u (s\ra u)^{-1}}\ ,
\end{eqnarray*}
again as required.
\endproof

This means that $kX\cobicross k(Y)$ inherits many of the nice properties of
$D(H)$. In particular, it has a quasitriangular structure $\CR$ and  associated
elements $Q=\CR_{21}\CR$ (the `quantum inverse Killing form') and $\cu=\sum
(S\CR\ut)\CR\uo$ (the element which implements the square of the antipode) in
Drinfeld's general theory of quasitriangular Hopf algebras\cite{Dri}.

\begin{corol} The bicrossproduct $kX\cobicross k(Y)$ is quasitriangular, with
\[ \CR=\sum_{u,s,v,t} v^{-1}\tens\delta_{us}\tens s^{-1}\tens \delta_{(s\la
v)t}.\]
The `quantum inverse Killing form' is
\[ \CR_{21}\CR\ =\
\sum_{u,v\in G\ s,p\in M}
s^{-1}u^{-1} \tens \delta_{u(s\la v)(p\ra u)u^{-1}} \tens v^{-1}p^{-1}
 \tens \delta_{pusp^{-1} } \ .\]
and is non-degenerate as a bilinear functional $k(X)\bicross k Y$. The element
$\cu$ is the canonical
element $\cu=\sum_{x\in X}x\tens\delta_x$
in $k X\cobicross k(Y)$, and is central.
\end{corol}
\proof The computation is straightforward. Thus
\align{(\psi\tens\psi)(\CR)\equad &&= \sum \psi(\delta_s\tens u\tens
e\tens\delta_v)\tens\psi(\delta_t\tens e\tens s\tens\delta_u)\\
&&=\sum(s\la u)^{-1}\tens\delta_{vs^{-1}}\tens s\tens \delta_{u(s\ra
u)^{-1}t^{-1}s}\\
&&=\sum (s^{-1}\la v)^{-1}\tens\delta_{us}\tens s^{-1}\tens\delta_{vt}}
which yields the formula shown on a change of variables $v$ to $s\la v$. We
then compute $\CR_{21}\CR$ using the product in $kX\cobicross k(Y)$.
\begin{eqnarray*}
\CR_{21}\CR & = & \Big(
\sum_{u,s,v,t}  s^{-1}\tens \delta_{(s\la
v)t} \tens v^{-1}\tens\delta_{us}\Big)\Big(
\sum_{u',s',v',t'} v'^{-1}\tens\delta_{u's'}\tens s'^{-1}\tens \delta_{(s'\la
v')t'}\Big) \\
& = & \sum_{u,s,v,t,u',s',v',t'}
 \Big(s^{-1}\tens \delta_{(s\la
v)t}\Big)\Big(v'^{-1}\tens\delta_{u's'}\Big)\tens\Big(
v^{-1}\tens\delta_{us}\Big)\Big(s'^{-1}\tens \delta_{(s'\la
v')t'}\Big) \\
& = & \sum_{u,s,v,t,u',s',v',t'}
\delta_{(s\la v)t, v'^{-1}\tilde\la u's'}\
\delta_{us, s'^{-1}\tilde\la (s'\la v')t'}\
s^{-1} v'^{-1}\tens\delta_{u's'} \tens
v^{-1} s'^{-1}\tens \delta_{(s'\la
v')t'}
\\
& = & \sum_{u,s,v,t,u',s',v',t'}
\delta_{(s\la v)t, v'^{-1} u'(s'\la v')(s'\ra v')}\
\delta_{us, v' (s'\ra v')^{-1}t's'}\
s^{-1} v'^{-1}\tens\delta_{u's'} \tens
v^{-1} s'^{-1}\tens \delta_{(s'\la
v')t'}
\end{eqnarray*}
{}From the $\delta$-functions here we can read off
$s\la v= v'^{-1} u'(s'\la v')$,
$t= s'\ra v'$,
$u= v' $, and
$s=  (s'\ra v')^{-1}t's'$. If we rewrite these as
$v'=u$, $s'=t\ra u^{-1}$, $u'=u(s\la v)(t\la u^{-1})$ and
$t'=ts(t\ra u^{-1})^{-1}$, then
on substituting $p=t\ra u^{-1}$:
\[
\CR_{21}\CR\ =\  \sum_{u,v\in G\ s,p\in M}
s^{-1}u^{-1} \tens \delta_{u(s\la v)(p\ra u)u^{-1}} \tens v^{-1}p^{-1}
 \tens \delta_{pusp^{-1} } \ .
\]
Nondegeneracy of $\CR_{21}\CR$ as a linear map $D(H)^*\to D(H)$ is the
so-called
factorisability property holding for any quantum double\cite{Ma:book}. Hence it
carries over in our case to a linear isomorphism $k(X)\bicross kY\to
kX\cobicross k(Y)$ or, equivalently, to a nondegenerate bilinear functional on
$k(X)\bicross k Y$.

Finally,
\align{\psi(\cu)\equad &&=\sum \psi(\delta_s\tens u\tens (s\ra
u)^{-1}\tens\delta_{(s\la u)^{-1}})\\
&&=\sum(s\la u)^{-1}(s\ra u)^{-1}\tens \delta_{(s\la u)^{-1}((s\ra
u)^{-1}\ra(s\la u)^{-1})^{-1}s^{-1}(s\ra u)^{-1}}\\
&&=\sum (s\la u)^{-1}(s\ra u)^{-1} \tens\delta_{(s\la u)^{-1}(s\ra
u)^{-1}}=\sum us\tens\delta_{us}}
on a change of summation variables in which  $s\la u$ is replaced by $s^{-1}$
and $s\la u$ by $u^{-1}$.
Finally, the element $\cu$ in any quasitriangular Hopf algebra implements the
square of the antipode. But for any bicrossproduct, the antipode is involutive,
hence $\cu$ here is central. \endproof

We also consider the $*$-structure in the case where the ground field is
equipped with an involution.

\begin{corol}\label{antireal} $k X\cobicross k(Y)$ is antireal-quasitriangular
in the sense $(*\tens *)(\CR)=\CR^{-1}$.
\end{corol}
\proof This is known for the quantum double $D(H)$ of any Hopf
$*$-algebra\cite{Ma:mec}, and hence follows from Theorem~\ref{Dbic}; here we
provide a direct proof for our particular case. Computing using the
bicrossproduct $*$-structure, we have
\align{(*\tens *)(\CR)\equad&&=\sum_{u,s,v,t}v\tens
\delta_{v^{-1}\tilde\la(us)}\tens s\tens\delta_{s^{-1}\tilde\la((s\la
v)t)}=\sum_{u,s,v,t}v\tens\delta_{v^{-1}u(s\la v)(s\ra v)}\tens s\tens
\delta_{v(s\ra v)^{-1}ts}\\
&&=\sum_{u,s,v,t}v\tens\delta_{u(s\ra v)}\tens s\tens \delta_{vt}}
after a change of $u,t$ variables in the last step.
Meanwhile,
\align{\CR^{-1}\equad &&=(S\tens\id)(\CR)=\sum_{u,s,v,t}(v^{-1}\tilde\ra
us)^{-1}\tens\delta_{(v^{-1}\tilde\la us)^{-1}}\tens s^{-1}\tens\delta_{(s\la
v)t}\\
&&=\sum_{u,s,v,t} s\la v\tens\delta_{(v^{-1}u(s\la v)(s\ra v))^{-1}}\tens
s^{-1}\tens\delta_{(s\la v)t}=\sum_{u,s,v,t} s\la v\tens\delta_{(s\la
v)^{-1}u^{-1}v(s\ra v)^{-1}}\tens s^{-1}\tens\delta_{(s\la v)t}}
using the actions in Propsition~2.1. Note that the inversion in $\delta_{(\
)^{-1}}$ is the inverse in $Y=G\times M^{\rm op}$. This gives the same as
$(*\tens *)(\CR)$ after changing variables to $v'=s\la v$, $s'=s^{-1}$,
$u'=(s\la v)^{-1}u^{-1}v$. Here $(s\ra v)^{-1}=s^{-1}\ra(s\la v)=s'\ra v'$.
\endproof

As an application, the finite-dimensional modules of any quasitriangular Hopf
algebra have a natural `quantum dimension' $\underline\dim$ defined as the
trace of $\cu$ in the representation. The modules of $kX\cobicross k(Y)$, as a
cross product algebra, are just the $Y$-graded $X$-modules $V$ such that
$|x\la v|=x\tilde\la |v|$ for all $v\in V$ homogeneous of degree $|\ |$.

\begin{propos} The quantum dimension of a general $kX\cobicross k(Y)$-module
$V$ is
\[  \underline\dim(V)\ =\ \sum_{y\in Y} {\rm trace}_{V_y}\,
\pi(y) \]
where $V_y$ is the subspace of degree $y$ and $\pi(y):V_y\to V_y$ is the
restriction to $V_y$ of the action of $y$ viewed as an element of $X$.
\end{propos}
\proof We write $V=\oplus_y V_y$ for our $Y$-graded $X$-module. The action of
$f\in k(Y)$ is $f(y)$ on $V_y$. A general element $x\in X$ acting on $V$ sends
$V_y\to V_{x\tilde\la y}$. Hence, in particular, $y$ viewed in $X$ sends
$V_y\to V_y$ as $y\tilde\la y=y$ from Proposition~2.1. This is the operator on
$V_y$ denoted $\pi(y)$. Let $\{e_a^{(y)}\}$ be a basis of $V_y$, with dual
basis $\{f^{a(y)}\}$. Then
\[ \trace(\cu)=\sum_{y\in Y,x\in X}\sum_{a} \<f^{a(y)}, (x\tens
\delta_x).e^{(y)}_a\>=
\sum_{y\in Y }\sum_{a} \<f^{a(y)}, y.e^{(y)}_a\>=\sum_{y\in
Y}\trace_{V_y}\pi(y).\]
\endproof

For example, we may take the natural representation in $k(Y)$ by left
multiplication of $k(Y)$ and the left action of $X$ induced by its action on
$Y$. This is the so-called Schroedinger representation of any cross product
algebra. The spaces $V_y$ are 1-dimensional with basis $\{\delta_y\}$ and
$\pi(y)\delta_y= \delta_y(y^{-1}\tilde\la(\ ))=\delta_{y\tilde\la
y}=\delta_{y}$ is the identity. So $\underline\dim(k(Y))=|Y|=\dim k(Y)$, where
$|Y|$ is the order of group $Y$. So for this representation the quantum
dimension is the usual dimension.

\begin{example}
We consider the factorisation of the group $S_3$ into a
subgroup of order 3 and a subgroup of order 2.
\end{example}
\proof Consider a factorisation of the group $S_3$ of permutations of 3
objects,
which we label 1,2 and 3. Let $G$ be the subgroup consisting of the
3-cycles and the identity, and let
the subgroup $M$ consist of the transposition $(1,2)$ and the identity. Then,
in the notation of this section, $X=S_3$, and $Y=GM^{op}$ is a cyclic group
of order 6.
The left action of $X$ on $Y$ is the adjoint action of the group $S_3$ on
the set
$S_3$, and the right action of $Y$ on $X$ is given by
\[
u\tilde\ra v=u\ ,\ u\tilde\ra v(1,2)=u^{-1}
\ ,\ u(1,2)\tilde\ra v=u(1,2)\ ,\ u(1,2)\tilde\ra v(1,2)=u^{-1}(1,2)\ ,
\]
where $u$ and $v$ are any 3-cycles or the identity. This leads to a
quasitriangular structure $\CR$
on $kX\cobicross k(Y)$, given by
\[
\CR\ =\ \sum_{u,v\in G,\ t\in M} v^{-1}\tens\delta_u\tens e\tens\delta_{vt}\ +\
\sum_{u,v\in G,\ t\in M} v^{-1}\tens\delta_{u(1,2)}\tens
(1,2)\tens\delta_{v^{-1}t}\ .
\]
What actually is the group  $Y\dcross X$ in this case? It is of order 36,
and a short
calculation will show that it has no center. The possibilities, read off
from a table of
groups, are $S_3\times S_3$ and $C_4\rcross (C_3\times C_3)$, where the $C_4
=\big\{\underline 0,\underline 1,\underline 2,\underline 3\big\}  $
 action is given by
$\underline 1\la (x,y)=(y,-x)$. The group $S_3\times S_3$ has 15 elements
of order 2,
and $C_4\rcross (C_3\times C_3)$ has 9 elements of order 2. A brief check
of the group $Y\dcross X$ shows that it has 15 elements of order 2, so it
is isomorphic to
$S_3\times S_3$. However the isomorphism does not seem to be obvious! \endproof

\section{Subfactorisation from an order-reversing isomorphism}

Let $\theta$ be an automorphism of $X$ which reverses its factors $GM$ (i.e.
$\theta(G)=M$ and $\theta(M)=G$). It is shown in \cite{BGM:fin} that $\theta$
induces an semi-skew automorphism of $D(H)$ (i.e. an algebra antiautomorphism
and
coalgebra automorphism), which we denote $\Theta$:
\eqn{DTheta}{\Theta(\delta_s\tens u\tens t\tens\delta_v)=\delta_{\theta(t\la
v)}\tens \theta(t\ra v)\tens \theta(s\la u)\tens\delta_{\theta(s\ra u)}.}
Via Theorem~\ref{Dbic}, we may view this as a semi-skew automorphism of
$kX\cobicross k(Y)$. When the ground field is equipped with an involution, we
may follow $\Theta$ by the star operation and obtain an antilinear Hopf algebra
automorphism $*\Theta$.

\begin{lemma} If $\theta$ is a factor-reversing automorphism of $X$ then the
induced antilinear automorphism of $k X\cobicross k(Y)$ is given by
 \[*\Theta(x\tens\delta_y)=\theta(x)\tens\delta_{\theta(y)^{-1}}\]
when $y$ is viewed in $X$ (and the inverse is also in $X$).
\end{lemma}
\proof We define $*\Theta$ via $\psi$ and (\ref{DTheta}). Thus,
\align{&&\nquad * \psi\big(\delta_{\theta(t\la
v)}\tens \theta(t\ra v)\tens \theta(s\la u)\tens\delta_{\theta(s\ra u)}\big)\\
&& =  *\Big(
\big(\theta(t\la v)\la\theta(t\ra v)\big)^{-1} \theta(s\la u) \tens
\delta_{ \theta(s\ra u) \big(\theta(s\la u)\ra\theta(s\ra u)\big)^{-1}
\theta(t\la v)^{-1}\theta(s\la u)}\Big) \\
&& =  *\Big(
\theta( t^{-1} (s\la u) )\tens \delta_{\theta\big((s\ra u) u^{-1}(t\la
v)^{-1}(s\la u)\big)}\big)  =
\theta(  (s\la u)^{-1}t )\tens \delta_{\theta\big(t^{-1}s
(t\ra v)^{-1}v^{-1}\big)}\\
&&=*\Theta\psi(\delta_s\tens u\tens t\tens\delta_v)=*\Theta\big((s\la u)^{-1}t
\tens \delta_{v(t\ra v)^{-1}s^{-1}t}\big).}
Comparing these expressions  gives the result for $*\Theta$ after changing
variables to general elements of $X,Y$.
 \endproof

We observe that $*\Theta$-invariant basis elements $x\tens\delta_y$
are characterised by the property that $\theta x=x$ and $\theta(y)=y^{-1}$
(computed in $X$).

\begin{propos} There is a subgroup $X^\theta$ of $X$ consisting of those
elements
$x$ for which $\theta x=x$, and a subset $Y^\theta$
 of $Y$ consisting of those elements
$y$ for which $\theta y=y^{-1}$ (inverse in $X$). The actions
$\tilde\la,\tilde\ra$ restrict to $X^\theta,Y^\theta$, forming a double cross
product group $Y^\theta\dcross X^\theta$ factorising into $Y^\theta,X^\theta$.
The corresponding bicrossproduct $kX^\theta\cobicross k(Y^\theta)$ Hopf algebra
has an isomorphic  coalgebra to that of $kM\cobicross k(G)$.
\end{propos}
\proof
The proof that the $\tilde\la$ action restricts is immediate. If we take
$x\in X^\theta$ and $y\in Y^\theta$, then $x \tilde\la y=xyx^{-1}$ (adjoint
action in the $X$ multiplication). If we apply $\theta$ to this, then
\[
\theta(x \tilde\la y)\ =\ \theta(xyx^{-1})\ =\
\theta(x)\theta(y)\theta(x^{-1})\ =\ xy^{-1}x^{-1}\ ,
\]
which in the inverse (in $X$) of $x \tilde\la y=xyx^{-1}$, so
$x \tilde\la y\in Y^\theta$.

The proof for the other action is rather more difficult.
It will be most convenient to find formulae for the elements of
 $X^\theta$ and $Y^\theta$ first.

If $y=vt\in Y^\theta$, then $\theta y=y^{-1}$ (inverse in $X$), so
we can substitute $\theta (y)=\theta (v)\theta (t)$ and
$y^{-1}=t^{-1}v^{-1}$ and use uniqueness of factorisation
to say that $\theta (v)=t^{-1}$. Then we can write
$y=Y(v)=v\theta(v)^{-1}$. Now we can write out
a simple formula for the multiplication in $Y^\theta$ as
\[
y.y'\ =\ Y(v).Y(v')\ =\ (v\theta(v)^{-1}).(v'\theta(v')^{-1})\ =\
vv'\theta(v')^{-1}\theta(v)^{-1}\ =\
vv'\theta(vv')^{-1}\ =\ Y(vv')\ .
\]
This shows that $Y^\theta$ is actually isomorphic to $G$.

If $x=us\in X^\theta$, then $\theta x=x$, so
we can substitute
\[\theta (x)\ =\ \theta (u)\theta (s)\ =\
 x\ =\ us\ =\ (s^{-1}\ra u^{-1})^{-1}\, (s^{-1}\la u^{-1})^{-1}\ ,
\]
 and use uniqueness of factorisation
to say that $\theta (s)=(s^{-1}\la u^{-1})^{-1}$.
Then $u^{-1}=s\la \theta(s)^{-1}$, so we can write
$x=X(s)=\big( s\la \theta(s)^{-1} \big)^{-1}s$. This shows that $X^\theta$ is
bijective as a set with $M$.

Now consider the right action,
 \begin{eqnarray*}
 X(s) \tilde\ra Y(v)& =&
\big( s\la \theta(s)^{-1} \big)^{-1}s \tilde\ra
v\theta(v)^{-1}\ =\ \big((s\ra v)
\theta(v)^{-1} s^{-1}\la(s\la \theta(s)^{-1})\big)^{-1}\ (s\ra v) \\
& = & \big((s\ra v)
\theta(v)^{-1}\la \theta(s)^{-1}\big)^{-1}\ (s\ra v)  \\
& = & \big((s\ra v) \la \theta(s\ra v)^{-1}
\big)^{-1}\ (s\ra v)\ =\ X(s\ra v)\ .
\end{eqnarray*}
In particular this shows that the result is in $X^\theta$.

We therefore have a bicrossproduct Hopf algebra $k X^\theta \cobicross
k(Y^\theta)$. Its coalgebra is determined by the action of $Y^\theta$ on
$X^\theta$ and the group structure of $Y^\theta$, hence we see that it is
isomorphic via the maps $X,Y$ to the coalgebra of $kM\cobicross k(G)$.\endproof

Also, by construction, we see that if we equip $k$ with the trivial involution
then
\[k X^\theta\cobicross k(Y^\theta)\subseteq (k X\cobicross k(Y))^{*\Theta}\]
as algebras. The right hand side denotes the fixed point subalgebra under the
algebra automorphism $*\Theta$. The inclusion is clear from Lemma~5.1
(????) and the
inclusion $X^\theta\subset X$ as subgroups and $k(Y^\theta)\subset k(Y)$
(extension by zero) as $X^\theta$-module algebras.
 That $*\Theta$ is also a
coalgebra automorphism tells us further that the coproduct of $k X\cobicross
k(Y)$ applied to elements of the fixed subalgebra yields elements invariant
under $*\Theta\tens *\Theta$. It is natural to ask to what extent the
quasitriangular structure of $k X\cobicross k(Y)$ is likewise invariant.

\begin{propos} The quasitriangular structure of $kX\cobicross k(Y)$ obeys
 \[ (\Theta\tens\Theta)(\CR)=\CR_{21}\]
When the field has an involution, we have $(*\Theta\tens
*\Theta)(\CR)=\CR_{21}^{-1}$. Moreover, if $\theta^2=\id$ then  $\Theta^2=\id$
and $(*\Theta)^2=\id$.
\end{propos}
\proof It is easier to do the first computations in $D(H)$. There, we have
\align{ (\Theta\tens\Theta)(\CR)\equad &&=\sum_{u,v,s,t}\Theta(\delta_s\tens
u\tens e\tens\delta_v)\tens\Theta(\delta_t\tens e\tens s\tens\delta_u)\\
&&=\sum_{u,v,s,t}(\delta_{\theta(e\la v)}\tens\theta(e\ra v)\tens\theta( s\la
u)\tens\delta_{\theta(s\la u)})\tens (\delta_{\theta(s\la u)}\tens \theta(s\ra
u)\tens \theta(t\la e)\tens\delta_{\theta(t\ra e)})\\
&&=\sum_{u,s}(1\tens\theta(s\la u)\tens\delta_{\theta(s\ra u)})\tens
(\delta_{\theta(s\la u)}\tens\theta(s\ra u)\tens 1)}
where the sum over $v,t$ are replaced by sums over $t'=\theta(v),v'=\theta(t)$
and give the unit elements of $k(M)$ and $k(G)$ respectively. Then we change
variables from $u,s$ to $s'=\theta(s\la u),u'=\theta(s\ra u)$, to recognise
$\CR_{21}$ in $D(H)$. Hence the same result applies for $kX\cobicross k(Y)$.
This combines with Corollary~\ref{antireal} to obtain the corresponding
property for $*\Theta$. Also, the automorphism
$\Theta$ in (\ref{DTheta}) clearly obeys
\[  \Theta^2(\delta_s\tens u\tens
t\tens\delta_v)=\delta_{\theta^2(s)}
\tens\theta^2(u)\tens\theta^2(t)\tens\delta_{\theta^2(v)}\]
and hence $\Theta^2=\id$ when $\theta^2=\id$. The same feature for $*\Theta$ is
immediate from Lemma~3.1. \endproof

Hence $\CR$ is not invariant in the usual sense (unless $G,M$ are trivial), due
to the nondegeneracy of $\CR_{21}\CR$ in  Corollary~2.3.
Rather, one should note that for any quasitriangular Hopf algebra,
$\CR_{21}^{-1}$ defines a second `conjugate' quasitriangular structure. It
corresponds in topological applications to reversing braid crossings; we see
that our quasitriangular structure is invariant up to conjugation in this
sense.
Although $kX^\theta\cobicross k(Y^\theta)$ does not in general inherit a
quasitriangular
structure from $\CR$,  its inclusion as fixed point subalgebra in a
quasitriangular
Hopf algebra equipped with such an automorphism might be a useful substitute.
Of course, it may still happen that $kX^\theta\cobicross k(Y^\theta)$ is
quasitriangular
for some other reason, which is the case in the first example below.

Another natural question, in view of Proposition~3.2, is whether $k
X^\theta\cobicross k(Y^\theta)$ is in fact isomorphic as a Hopf algebra to our
original bicrossproduct $kM\cobicross k(G)$. The following example also
demonstrates
that it is not necessarily isomorphic to it.

\begin{example}
We consider the example in [BGM] of the doublecrossproduct of two
cyclic groups of order 6 ($C_6$) which gives the product of two symmetric
groups $S_3\times S_3$. In this case
$kX^\theta\cobicross k(Y^\theta)$ is isomorphic to $kS_3\cobicross k(S_3)$ in
Example~2.6 and
hence quasitriangular.
\end{example}
\proof Consider the group
$X=S_3\times S_3$ as the permutations of 6 objects labelled 1 to 6,
where the
first factor leaves the last 3 objects unchanged, and the second
factor leaves the first 3 objects unchanged. We take $G$ to be the
cyclic group
of order 6 generated by the permutation $1_G=(123)(45)$, and
$M$ to be the cyclic group
of order 6 generated by the permutation $1_M=(12)(456)$.
Our convention is that
permutations act on objects on their right, for example $1_G$
applied to $1$
gives $2$. The intersection of $G$ and $M$ is just the identity
permutation, and
counting elements shows that $GM=MG=S_3\times S_3$. We write
each cyclic group
additively, for example $G=\{0_G,1_G,2_G,3_G,4_G,5_G\}$.
The action of the
element $1_M$ on $G$ is seen to be given by the permutation
$(1_G,5_G)(2_G,4_G)$, and that of  $1_G$ on $M$ is given by the
permutation
$(1_M,5_M)(2_M,4_M)$.

The factor reversing automorphism $\theta$ can be taken to be
conjugation by the permutation $(1,4)(2,5)(3,6)$. Then if we split the
elements of $X$ into $S_3\times S_3$, we see that
the elements of $X^\theta$ are of the form $\sigma\times\sigma$,
for $\sigma\in S_3$, and that the elements of $Y^\theta$
are of the form $\sigma\times\sigma^{-1}$. Then $X^\theta$ is isomorphic to
$S_3$, and the action of $X^\theta$ on $Y^\theta$ is the
adjoint action of the group $S_3$ on the set $S_3$.

This is enough information to show that, despite having the same
dimension, that
$kX^\theta\cobicross k(Y^\theta)$ is not isomorphic
to $kG\cobicross k(M)$ or $k(G)\bicross kM$.
All these can, as algebras
over $\Bbb C$, be decomposed into a  direct sum of matrix algebras.
According to \cite{BGM:fin}, the algebras $kG\cobicross k(M)$ or $k(G)\bicross
kM$
have at most $2\times 2$ matrices in their decompositions. However
the existence of an orbit of size 3 in the
adjoint action of the group $S_3$ on the set $S_3$ shows that
at least one $3\times 3$ matrix occurs in the decomposition
of $kX^\theta\cobicross k(Y^\theta)$.

So what is $kX^\theta\cobicross k(Y^\theta)$? In fact we have already met
it, it is
the bicrossproduct given by the product of $S_3$ and a cyclic group of
order 6 in example 2.6. The explicit correspondence is given by
deleting the points 4, 5 and 6 from the example here.
Then we have the maps $\sigma\times\sigma\mapsto \sigma$ for
$X^\theta$, and $\sigma\times\sigma^{-1}\mapsto \sigma$ for
$Y^\theta$, corresponding to the subsets of $S_3$
used in in example 2.6.
But from the previous example we see that $kX^\theta\cobicross k(Y^\theta)$
is actually the double of a bicrossproduct obtained from a group of order 2
and a group of
order 3. We deduce that as a result $kX^\theta\cobicross k(Y^\theta)$
is actually quasitriangular. However there is no obvious relation between this
quasitriangular structure and the standard structure on $kX\cobicross k(Y)$.
\endproof

\begin{example} The upper-lower triangular example, in which
$kX^\theta\cobicross k(Y^\theta)$ is isomorphic to $kM\cobicross k(G)$.
\end{example}
\proof Let $G$ be the group $T_+$ of upper triangular $n\times n$ matrices with
ones on the diagonal, with the usual matrix
multiplication. Also let $M$ be the
group $T_-$ of lower triangular $n\times n$ matrices with
ones on the diagonal.
As in \cite{Ma:book} we define actions, using $\theta(a)=(a^{T})^{-1}$,
\[
s\la u\ =\ 1\ +\ \theta(s)(u-1)\quad,\quad
s\ra u\ =\ 1 + (s-1)\theta(u)\ .
\]
Using these actions we find
  \[
X(t).X(s)\ =\
 X\big( (t\ra(s\la \theta(s)^{-1})^{-1})s\big)=\ X(2ts-s-t+1)\ .
\]
Hence there is a group isomorphism $X^\theta\isom T_-$
defined by $X(t)\mapsto 2t-I$. We call the group isomorphism
 $\alpha:T_-\to X^\theta$, where $\alpha(s)=X((s+1)/2)$. Then
we calculate
\[
\alpha(s)\tilde\ra Y(v)=X((s+1)/2)\tilde\ra Y(v)=X({s+1\over 2}\ra v)=
X\big(1+({s+1\over 2}-1)\theta(v)\big)=X\big(1+{s-1\over 2}\theta(v)\big)=
\alpha(s\ra v)\ .
\]
As the map $v\to Y(v)$ gives a group isomorphism from $T_+$ to $Y^\theta$
we might think that
we are on our way to showing that the $X^\theta Y^\theta $
doublecross product is isomorphic to the original
$T_+T_-$ doublecross product. To check this, we must perform the
calculation of the left action:
\[
X(s)\tilde\la Y(v)=Y\big(1+(2-\theta(s))^{-1}\theta(s)(v-1)\big)=
Y\big(1+(2\theta(s)^{-1}-1)^{-1}(v-1)\big)
\]
Now we use the fact that $\theta(s)^{-1}=s^T$, and calculate
\[
\alpha(t)\tilde\la Y(v)=X((t+1)/2)\tilde\la Y(v)=
Y\big(1+(t^T)^{-1}(v-1)\big)=Y(t\la v)\ ,
\]
as required to prove the isomorphism of the doublecross products. \endproof

\section{A *-representation of $D(H)$ on a
Hilbert space}
 In this section we provide a
Hilbert space representation of $D(H)$ which is
one of the motivations behind Theorem~\ref{Dbic}.  Recall that it was shown in
\cite{BGM:fin} that representations of $D(H)$ are $G-M$-bigraded bicrossed
$G-M$-bimodules. We shall use $|w|$ for the $G$-grading and $\big<w\big>$
for the $M$-grading
of a homogeneous element $w$ of the representation.

\begin{propos} There is a representation of $D(H)$ on a vector space  $E$
with
basis $\{\eta_{s,u}|s\in M, u\in G\}$, with gradings
\[
\big| \eta_{s,u} \big|\ =\ u \quad,\quad \big< \eta_{s,u} \big>
\ =\ s\
\]
and the group actions
\[
t\la  \eta_{s,u} \ =\ \eta_{ts(t\ra u)^{-1},t\la u}
 \quad,\quad
\eta_{s,u} \ra v \ =\ \eta_{s\ra v,(s\la v)^{-1}uv}\ .
\]
The corresponding action of $D(H)$ is
\[
(\delta_s \tens u\tens t\tens\delta_v)\la \eta_{r,w}\ =\
\delta_{v,w}\ \delta_{t^{-1}s(t\ra v),r}\ \eta_{
s\ra u,(s\la u)^{-1}(t\la w)u  }\ .
\]
\end{propos}
\proof The definition of the group actions
is made precisely so that the matching
conditions in \cite{BGM:fin} are true.
The corresponding actions of the Hopf algebras $H^{*}$ and $H$  are
\[
\big(t\tens\delta_v\big)\la \eta_{s,u} \ =\
\delta_{v,u}\ t\la  \eta_{s,u}
 \quad,\quad
\eta_{s,u}\ra \big(\delta_t\tens v\big) \ =\
\delta_{s,t}\ \eta_{s,u} \ra v\ ,
\]
and the formula
$(a\tens h)\la w=(h\la w)\ra a$  gives the action of $D(H)$ as
\[
\big(I\tens t\tens\delta_v \big)
\la \eta_{s,u} \ =\
\delta_{v,u}\
\eta_{ts(t\ra u)^{-1},t\la u}
 \quad,\quad
\big(\delta_t\tens v\tens I \big)
\la \eta_{s,u} \ =\
\delta_{t,s}\
\eta_{s\ra v,(s\la v)^{-1}uv}\ ,
\]
which gives the formula stated. \endproof

As far as the original group doublecross product is concerned, the
$E$ representation is more symmetric than the standard\cite{Ma:book}
`Schroedinger'
representation of $D(H)$ in $H$, as we do not have to decide to take the
group algebra of one factor
subgroup, and the function algebra of the other.
The $E$ representation
motivates the isomorphism in Theorem~\ref{Dbic} in the following manner:
If we rewrite $\rho_{wr^{-1}}=\eta_{r,w}$, then the action above
gives
\[
(\delta_s \tens u\tens t\tens\delta_v)\la \rho_{wr^{-1}}\ =\
\delta_{v(t\ra v)^{-1}s^{-1}t,wr^{-1}}\
\rho_{(s\la u)^{-1}t(wr^{-1})t^{-1}(s\la u)}\ .
\]
Compare this to the action of a single bicrossproduct
$kX\cobicross k(Y)$, with $Y$-grading $\|\ \|$ on homogeneous elements:
\[
(x\tens\delta_y)\la (\ ) =\ \delta_{y,\|w\|}\ x\la (\ ).
\]
We see that the formulae agree if we set $x=(s\la u)^{-1}t$,
 $y=v(t\ra v)^{-1}s^{-1}t$,
$\|\rho_{wr^{-1}}\|=wr^{-1}$,
and use the adjoint action of $X$, $x\la
\rho_{wr^{-1}}=\rho_{x(wr^{-1})x^{-1}}$.
 This suggests trying a formula of the
type used for $\psi$ in Section~2.

\begin{propos} Over $\C$, there is an inner product $\big(\eta,\zeta\big)$ on
$E$ (conjugate-linear in the first variable
 and linear on the second), defined by
\[
\big(\eta_{s,u},\eta_{r,w}\big)\ =\
\delta_{s,r}\ \delta_{u,w}\ .
\]
With respect to this, $D(H)$ with its natural $*$-structure is represented
as a $*$-algebra.
\end{propos}
\proof  The inner product given is non-degenerate (in fact the $\eta_{s,u}$
form an orthonormal basis).
 We can then check the condition that
\[
\big(\alpha\la \eta,\zeta\big)\ =\ \big(\eta,\alpha^*\la \zeta\big)
\]
for any $\alpha\in D(H)$. We shall only prove this in the case
$\alpha=I\tens t\tens\delta_v$, and leave the other case to the
reader. First we calculate
\[
\big(\alpha\la \eta_{su},\eta_{s'u'}\big)\ =\
\delta_{v,u}\ \delta_{ts(t\ra u)^{-1},s'}\ \delta_{t\la u,u'}\ .
\]
Now $\alpha^*=I\tens t^{-1}\tens\delta_{t\la v}$, and
\[
\big(\eta_{su},\alpha^*\la \eta_{s'u'}\big)\ =\
\delta_{t\la v,u'}\ \delta_{s,t^{-1}s'(t^{-1}\ra u')^{-1}}\
\delta_{u,t^{-1}\la u'}\
, \]
which is the same condition. \endproof

\begin{propos}  $E$  has a coalgebra structure
\[
\Delta_E\big(\eta_{s,u}\big)\ =\
\sum_{ab=s,\ wz=u} \eta_{a,w}\tens \eta_{b,z}\ =\
\sum \eta_{(1)}\tens  \eta_{(2)} \quad , \quad \eps\big(\eta_{s,u}\big)\ =\
\delta_{s,e}\ \delta_{u,e}\ ,
\]
and  becomes a left module coalgebra under the action of $D(H)$
\end{propos}
\proof  To show that $E$ is a coalgebra, we first show that it is
coassociative, which is
\[
(I\tens\Delta_E)\Delta_E\big(\eta_{s,u}\big)\ =\
\sum_{abc=s,\ wzv=u} \eta_{a,w}\tens \eta_{b,z}\tens \eta_{c,v}\ =\
(\Delta_E\tens I)\Delta_E\big(\eta_{s,u}\big)\ .
\]
Next we show that $\eps$ is a counit;
\[
(I\tens\eps)\Delta_E\big(\eta_{s,u}\big)\ =\
\sum_{ab=s,\ wz=u} \eta_{a,w}\ \delta_{e,z}\delta_{b,e}\ =\
\eta_{s,u}\ =\ (\eps\tens I)\Delta_E\big(\eta_{s,u}\big)\ .
\]
For the module coalgebra condition, first we have to prove
\[
\Delta_E(\alpha\la \eta)\ =\ \sum \alpha_{(1)}\la \eta_{(1)} \tens
\alpha_{(2)}\la \eta_{(2)}\ =\ \Delta(\alpha) \la \Delta_E(\eta)\ .
\]
Again we shall only prove this in the case
$\alpha=I\tens t\tens\delta_v$, and leave the other case to the
reader. We set $\eta=\eta_{s,u}$.  The coproduct of $D(H)$ is
\[
\Delta(\alpha)\ =\ \sum
\alpha_{(1)}\tens \alpha_{(2)}\ =\ \sum_{xy=v}
\big( I\tens t\tens\delta_x \big)
\tens
\big( I\tens t\ra x\tens\delta_y \big) \ ,
\]
and so we obtain
\[
\sum \alpha_{(1)}\la \eta_{(1)} \tens
\alpha_{(2)}\la \eta_{(2)}\ =\ \delta_{v,u}\ \sum_{xy=u\ ,\ ab=s}
\eta_{ta(t\ra x)^{-1},t\la x}
\tens
\eta_{(t\ra x)b(t\ra u)^{-1},(t\ra x)\la y}\ .
\]
Now we use $\alpha\la \eta=\delta_{v,u}\
\eta_{ts(t\ra u)^{-1},t\la u}$ and calculate
\[
\Delta_E(\alpha\la \eta)\ =\ \delta_{v,u} \sum_{x'y'=t\la u\ ,\ a'b'=ts(t\ra
u)^{-1}}
\eta_{a',x'}\tens \eta_{b',y'}\ .
\]
If we now use the correspondences $a'=ta(t\ra x)^{-1}$, $b'=(t\ra x)b(t\ra
u)^{-1}$,
$x'=t\la x$ and $y'=(t\ra x)\la y$ we see that the formulae for
$\alpha_{(1)}\la \eta_{(1)} \tens
\alpha_{(2)}\la \eta_{(2)}$ and $(\alpha\la \eta)_{(1)}\tens (\alpha\la
\eta)_{(2)}$ agree.

Lastly we must prove that
$
\eps(\alpha\la \eta) = \eps(\alpha)\eps(\eta)
$.
Using $\alpha=\delta_s \tens u\tens t\tens\delta_v$ and $\eta=\eta_{r,w}$,
we have
\begin{eqnarray*}
\eps\big(\alpha\la \eta\big)& =& \delta_{v,w}\ \delta_{t^{-1}s(t\ra v),r}\
\eps\big(\eta_{
s\ra u,(s\la u)^{-1}(t\la w)u  }\big)\ =\
\delta_{v,w}\ \delta_{t^{-1}s(t\ra v),r}\
\delta_{s\ra u,e}\ \delta_{(s\la u)^{-1}(t\la w)u,e}   \\
& =& \delta_{v,w}\ \delta_{t^{-1}s(t\ra v),r}\
\delta_{s,e}\ \delta_{(s\la u)^{-1}(t\la w)u,e} \ =\
\delta_{v,w}\ \delta_{t^{-1}(t\ra v),r}\
\delta_{s,e}\ \delta_{u^{-1}(t\la w)u,e}  \\
& =& \delta_{v,w}\ \delta_{t^{-1}(t\ra v),r}\
\delta_{s,e}\ \delta_{(t\la w),e} \ =\
 \delta_{v,w}\ \delta_{t^{-1}(t\ra v),r}\
\delta_{s,e}\ \delta_{w,e} \\
& =&  \delta_{v,e}\ \delta_{t^{-1}(t\ra v),r}\
\delta_{s,e}\ \delta_{w,e} \ =\  \delta_{v,e}\ \delta_{e,r}\
\delta_{s,e}\ \delta_{w,e} \\
& =&  \eps\big(\delta_s \tens u\tens t\tens\delta_v\big)\
\eps\big(\eta_{r,w}\big)
\ .
\end{eqnarray*}
\endproof

Now suppose that $\theta$ is an order reversing
isomorphism of the doublecross product group $X=GM$.
{}As previously mentioned, from \cite{BGM:fin} we have an
antilinear Hopf algebra
automorphism $*\Theta:D(H)\to D(H)$.

\begin{propos} There is an antilinear map $\hat\theta:E\to E$ defined by
\[
\hat\theta(\eta_{s,u})\ =\ \eta_{\theta(u),\theta(s)}\ ,
\]
which obeys
\[
\hat\theta\big(\alpha\la \eta\big)\ =\ (*\Theta \alpha)\la\hat\theta(\eta)\
\quad\forall \alpha\in D(H).
\]
\end{propos}
\proof
As usual we shall only prove this in the case
$\alpha=I\tens t\tens\delta_v$, and leave the other case to the
reader. We begin with
\[
*\Theta
\big(I\tens t\tens\delta_v\big)\ =\
*\big(\delta_{\theta(t\la v)} \tens \theta(t\ra v) \big) \tens I\ =\
\delta_{\theta(v)} \tens \theta(t\ra v)^{-1}\tens I\ ,
\]
where we have used the equation $\theta(t\la v)\ra \theta(t\ra v)
=\theta(v)$. Now
\[
\big(\delta_{\theta(v)} \tens \theta(t\ra v)^{-1}\tens I\big)
\la \eta_{\theta(u),\theta(s)}\ =\
\delta_{\theta(v),\theta(u)}\ \eta_{
\theta(u)\ra\theta(t\ra v)^{-1}\ ,\
\big( \theta(u)\la\theta(t\ra v)^{-1} \big)^{-1}
\theta(s)\theta(t\ra v)^{-1}  }\ .
\]
Since $\theta$ is 1-1, we can replace $\delta_{\theta(v),\theta(u)}$ by
$\delta_{v,u}$. Also we calculate
\[
\theta(u)\la\theta(t\ra v)^{-1}\ =\
\theta\big((t\ra v)\ra u^{-1}\big)^{-1}
\quad \text{and}\quad
\theta(u)\ra\theta(t\ra v)^{-1}\ =\
\theta\big((t\ra v)\la u^{-1}\big)^{-1}\ ,
\]
so the equation above becomes
\[
\big(\delta_{\theta(v)} \tens \theta(t\ra v)^{-1}\tens I\big)
\la \eta_{\theta(u),\theta(s)}\ =\
\delta_{v,u}\
 \eta_{
\theta\big((t\ra v)\la u^{-1}\big)^{-1}\ ,\
\theta\big((t\ra v)\ra u^{-1}\big)
\theta(s)\theta(t\ra v)^{-1}  }\ .
\]
Now using the $\delta_{v,u}$ to put $v=u$ in the equation;
\[
\big(\delta_{\theta(v)} \tens \theta(t\ra v)^{-1}\tens I\big)
\la \eta_{\theta(u),\theta(s)}\ =\
\delta_{v,u}\
 \eta_{
\theta(t\la u)\ ,\
\theta(t)
\theta(s)\theta(t\ra u)^{-1}  } \ =\
\delta_{v,u}\
 \eta_{
\theta(t\la u)\ ,\
\theta(ts(t\ra u)^{-1})  }
\ ,
\]
which is the formula for $\hat\theta\big(\alpha\la \eta_{s,u}\big)$
as required. \endproof

\begin{propos} The coproduct on $E$ and the inner product are related
 by the formula
\[
\big(\Delta_E \eta\, ,\, \Delta_E \zeta\big)\ =\ |X|\
\big( \eta\, ,\,  \zeta\big)\ ,
\]
where we use the tensor product inner product on $E\tens E$.
\end{propos}
\proof It is sufficient to prove this for basis elements,
\begin{eqnarray*}
\big(\Delta_E \eta_{s,u}\, ,\, \Delta_E \eta_{s',u'}\big)& = &
\sum_{ab=s,\ a'b'=s',\ wz=u,\ w'z'=u'}
\big( \eta_{a,w}\, ,\,  \eta_{a',w'}\big)\
\big( \eta_{b,z}\, ,\,  \eta_{b',z'}\big)  \\
& = &
\sum_{ab=s,\ a'b'=s',\ wz=u,\ w'z'=u'}
\delta_{a,a'}\  \delta_{w,w'}\  \delta_{b,b'}\  \delta_{z,z'} \\
& = &
\sum_{ab=s=s',\ wz=u=u'}
1 \ =\
\delta_{s,s'}\  \delta_{u,u'} \ |G|\ |M|
\ =\ |X|\
\big( \eta_{s,u}\, ,\,  \eta_{s',u'}\big)\ .
\end{eqnarray*}
\endproof

\section{First order bicovariant differential calculi on $H$}

In this section, we regard the Hopf algebra $A=H^*=k(M)\bicross k G$
associated to a group factorisation
$GM$ as a `coordinate ring' of some non-commutative geometric phase space. This
is the point of view introduced in \cite{Ma:phy}, where $H^*$ is an algebraic
model for the quantization of a particle on $M$ moving along orbits under the
action of $G$. Here we develop some of the `noncommutative geometry' associated
to this point of view.

First of all, on any algebra $A$, one may define a first-order   differential
calculus or `cotangent space' $\Omega^1$ in a standard way cf\cite{Con:non}

1. $\Omega^1$ is an $A$-bimodule.

2. $\extd:A\to \Omega^1$ is a linear map obeying the Leibniz rule
$\extd(ab)=(\extd a)b+a\extd b$.

3. $A\tens A\to \Omega^1$ given by $a\tens b\mapsto a\extd b$ is surjective.

When $A$ is a Hopf algebra, it is natural to add to this left and right
covariance (bicovariance) under $A$. Thus\cite{Wor:dif}

4. $\Omega^1$ is an $A$-bicomodule and the the bimodule structures and
$\extd$ are bicomodule maps. Here  $\Omega^1\tens A$ and $A\tens\Omega^1$
have the induced tensor product bicomodule structures, where $A$ is a
bicomodule under its coproduct.

Cf. \cite{Nic:bia} one knows that compatible bimodules and bicomodules (Hopf
bimodules) are of the form (say) $\Omega^1=V\tens A$ for some left crossed
$A$-module $V$. The latter in our finite-dimensional setting means nothing
more than left modules of the Drinfeld quantum double $D(A)$. A particular
module is $\ker\eps\subset A$, a restriction of the canonical `Schroedinger'
representation whereby $D(A)$ acts on $A$ (by left multiplication and the
coadjoint action, see\cite{Ma:book}). As observed in \cite{Wor:dif}, the
further conditions for $\Omega^1,\extd$ amount to requiring $V$ to be a
quotient of $\ker\eps$ as a quantum double module. Then
\eqn{extd}{ \extd a= (\pi\tens\id)(a\o\tens a\t-1\tens a),\quad\forall a\in A,}
where $\pi:\ker\eps\to V$ is the quotient map. The right (co)module structure
on $V\tens A$ is by right (co)multiplication in $A$, the left module
structure is the tensor product of $V$ and left (co)multiplication in $A$.

In the finite-dimensional setting which concerns us, one may equally well work
in the dual picture in terms of $H=k M\cobicross k(G)$ and $L=V^*$. By
definition cf\cite{Wor:dif}, a {\em bicovariant quantum tangent space} $L$ for
$H$ is a submodule of $\ker\eps\subset H$ under the
quantum double $D(H)$. Here $D(H)$ acts on $\ker\eps$ by
\[ h\la g=h\o g Sh\t,\quad a\la g=\<a,g\o\>g\t - \<a,g\>1,\quad \forall h\in
H,\quad g\in\ker\eps\subset H,\quad a\in H^*\]
as a projection to $\ker\eps$ of the Schroedinger representation. We say that a
quantum tangent space is irreducible if $L$ is irreducible as a quantum double
module. It corresponds to $\Omega^1$ having no bicovariant quotients. This dual
point of view has been used recently in \cite{Ma:cla}, where the irreducible
bicovariant quantum tangent spaces over a general quasitriangular Hopf algebra
have been classified under the assumptions of $\CR_{21}\CR$ non-degenerate and
a Peter-Weyl decomposition for the left regular representation. In the same
manner, we now classify the irreducible bicovariant quantum tangent spaces $L$
when $H$ is a bicrossproduct. The corresponding $\Omega^1$ will be given as
well.

The Schroedinger representation of $D(H)$ in $H$ has already been computed for
$H=k M\cobicross k(G)$ in \cite{BGM:fin}. Using the description of $D(H)$
modules as $M-G$ bicrossed bimodules (as recalled in Section~4), it
has\cite{BGM:fin}
\ceqn{Dact}{|t\tens\delta_v|=(t\la v)v^{-1},\quad \<t\tens\delta_v\>=t\\
s\la(t\tens\delta_v)=sts'^{-1}\tens \delta_{s'\la v},\quad (t\tens\delta_v)\ra
u=t\ra u\tens\delta_{u^{-1}v},}
where $s'=s\ra(t\la v)v^{-1}$.
We now use the isomorphism $D(H)\isom D(X)$ (the quantum double of the group
algebra of $X$) also in \cite{BGM:fin} to transfer to an action of $D(X)$.
Since $D(X)=k(X)\lcross kX$ {\em as an algebra}, this should make it easier to
decompose representations into irreducibles. One can also use our new
isomorphism in
Theorem~\ref{Dbic} to transfer to an action of $kX\cobicross k(Y)$, but this
appears to be less natural for the present purpose.

\begin{lemma} The action (\ref{Dact}) defines an action of  $k(X)\lcross
k X$ on $H$ given by
\[us\la (t\tens\delta_v)=(s''ts^{-1})\ra u^{-1}\tens\delta_{u(s\la v)},\quad
\delta_{us}\la(t\tens\delta_v)=\delta_{u,v((t\ra v)^{-1}\la
v^{-1})}\delta_{s\ra v,(t\ra v)^{-1}}t\tens\delta_v,\]
where $s''=s\ra v(t\la v)^{-1}$.
\end{lemma}
\proof  According to the general results in \cite{BGM:fin}, the corresponding
$X$-grading $||\ ||$ and $X$-action can be computed from the $M-G$-bicrossed
bimodule structure as $||w ||=\<w \>^{-1}| w |$ and $us\la w=((s\ra|w|^{-1})\la
w)\ra u^{-1}$ for all $w$ in the module. Hence,
\align{ ||t\tens\delta_v||\equad&&=t^{-1}(t\la v)v^{-1}=v((t^{-1}\ra(t\la
v))\la v^{-1}) (t^{-1}\ra((t\la v)v^{-1}))\\
&&=v((t\ra v)^{-1}\la v^{-1})((t\ra v)^{-1}\ra v^{-1})=v(t\ra v)^{-1}v^{-1}}
and
\cmath{us\la(t\tens\delta_v)=(s''\la (t\tens\delta_v))\ra
u^{-1}=(s''ts^{-1}\tens\delta_{s\la v})\ra u^{-1}=(s''ts^{-1})\ra
u^{-1}\tens\delta_{u(s\la v)}.}
\endproof

Motivated by the form of $||t\tens\delta_v||$ in the proof of the preceding
lemma, we chose new bases for the vector space on which $D(X)$ acts.

\begin{lemma} Let $\phi_{vt}=t^{-1}\ra v^{-1}\tens\delta_v$. Here
$\{\phi_{vt}\}$ is a basis of $H$ labelled by $vt\in X$. Then the action
in Lemma~5.1 is
\[ us\la \phi_{vt}=\phi_{usvt(s\ra v)^{-1}},\quad
\delta_{us}\la \phi_{vt}=\delta_{us,vtv^{-1}}\phi_{vt}\]
\end{lemma}
\proof Here $||\phi_{vt}||=vtv^{-1}=v(t\la v^{-1})(t\ra v^{-1})$ gives the
action of $\delta_{us}$ by evaluation against the degree. This can be written
more explicitly as $\delta_{us}\la\phi_{vt}=\delta_{u,v(t\la
v^{-1})}\delta_{s,t\ra v^{-1}}\phi_{vt}$ and is thereby equivalent to the
action in Lemma~5.1. Moreover,
\align{s''(t^{-1}\ra v^{-1})s^{-1}\equad&&=
(s\ra v((t^{-1}\ra v^{-1})\la v)^{-1})(t^{-1}\ra v^{-1})s^{-1}=(s\ra v
(t^{-1}\la v^{-1}))(t^{-1}\ra v^{-1})s^{-1}\\
&&=(((s\ra v)t^{-1})\ra v^{-1})s^{-1}=((s\ra v)t^{-1}(s\ra v)^{-1})\ra (s\la
v)^{-1}=((s\ra v)t(s\ra v)^{-1})^{-1}\ra (s\la v)^{-1}.}
Note that $(s\ra v)^{-1}\la (s\la
v)^{-1}=v^{-1}$ and $(s\ra v)^{-1}\ra (s\la v)^{-1}=s^{-1}$ for any matched
pair of groups. Then
\align{us\la\phi_{vt}\equad&&=us\la(t^{-1}\ra
v^{-1}\tens\delta_v)=(s''(t^{-1}\ra v^{-1})s^{-1})\ra u^{-1}\tens
\delta_{u(s\la v)}\\
&&=((s\ra v)t(s\ra v)^{-1})^{-1}\ra(u(s\la v))^{-1}\tens\delta_{u(s\la
v)}=\phi_{u(s\la v)(s\ra v)t(s\ra v)^{-1}}=\phi_{usvt(s\ra v)^{-1}}.}
\endproof

Our task is to
decompose $\ker\eps\subset H$ into irreducibles under this action of
$k(X)\lcross kX$.
We begin by decomposing the action in Lemma~5.2 into irreducibles and
afterwards projecting to $\ker\eps$. Note that Lemma~5.2 tells us that when we
identify $H\isom k X$ as linear spaces by the above basis, the action of $X$ is
the linear extension of a certain group action of $X$ on itself.

\begin{propos} Let $X$ act on itself by the action $us\tilde\la vt=us vt (s\ra
v)^{-1}$ as in Lemma~5.2. Let $||vt||\equiv ||\phi_{vt}||=vtv^{-1}$ as an
$X$-valued function on $X$.

(i) Let $\sim$ denote the equivalence relation on $X$ defined by $vt\sim us$
if and only if $||us||=||vt||$. Then $\tilde\la$ descends to an action of
$X$ on the
quotient space $X/\!\! \sim$.

(ii) Let  $\Xi_{[vt]}\subseteq X$ denote the isotropy subgroup of an
equivalence class $[vt]\in X/\!\! \sim$. Then
\[ \Xi_{[vt]}=\{us\in X|\ us ||vt|| =||vt||us\},\]
the centraliser of $||vt||$ in $X$.
\end{propos}
\proof (i) may be verified directly. However, it follows from Lemma~5.2 since
an action of $k(X)\lcross kX$ (where $X$ acts on $X$ by the adjoint action
in the semidirect product) requires
$||us\la\phi_{vt}||=us||\phi_{vt}||(us)^{-1}$. In terms of the group $X$, this
is $||us\tilde\la vt||=us||vt||(us)^{-1}$. This also implies (ii) since the
group $\Xi_{[vt]}$ consists of $us\in X$ such that $us\tilde\la vt\sim vt$,
i.e.
such that $||us\tilde\la vt||=||vt||$. \endproof

 We denote by $\CO_{[vt]}$ the orbit containing the point $[vt]$ in
$X/\!\!\sim$.

\begin{example} We may restrict attention to orbits of the form $\CO_{[s]}$,
where $s\in M$. Then the elements of the
equivalence class $[s]$ may be identified with the subset of $M$ fixed under
the action of $s$, $[s]=\big\{ sv|\ s\la v=v\big\}$.
The stabiliser $\Xi_{[s]}$ consists of all elements of $X$ which
commute with $s$.
The action of $\Xi_{[s]}$ on elements of the equivalence class $[s]$ is
given by $ut\tilde\la sv=
su(t\la v)$. In the particular case where $s=e$, we get $[e]=G$ and
$\Xi_{[e]}=X$.
\end{example}
\proof This may seem to be a rather specialised example, but in fact any orbit
$\CO_{[us]}$ in
$X/\!\!\sim$ contains a point of the given form, since $[s]\in \CO_{[us]}$. We
compute,
\[
[s]=\big\{ vt|\ vtv^{-1}=s\big\}=\big\{ vt|\ v(t\la v^{-1})=e\ ,\ t\ra
v^{-1}=s\big\}=
\big\{ v(s\ra v)|\ s\la v=v\big\}\ .
\]
This can be simplified if we note that if $s\la v=v$, then $v(s\ra v)=(s\la
v)(s\ra v)=sv$, giving the result stated. The action of $\Xi_{[s]}$ is computed
as
\[
ut\tilde\la sv=
ut\tilde\la v(s\ra v)=utv(s\ra v)(t\ra v)^{-1}=utsv(t\ra v)^{-1}=sutv(t\ra
v)^{-1}=
su(t\la v)\
\]
as stated. \endproof

For $\chi\in X/\!\!\sim$, define
 $\CS_\chi=k p^{-1}(\chi)\subset kX$, where $p$ is the canonical projection to
$X/\sim$. Here $\CS_\chi$ is the linear span of the elements of $\chi$ viewed
as a leaf in $X$, and is a linear $\Xi_\chi$ representation since, by
definition, the action of $\Xi_\chi$
sends $p^{-1}(\chi)$ to itself.

\begin{propos} Let $\CO$ be an orbit in $X/\!\!\sim$ under the action of
$X$. Then
$\CM_{\CO}=\oplus_{\chi\in\CO}\CS_\chi\subset kX$ is a subrepresentation under
the action of $k(X)\lcross kX$ in Lemma~5.2. Morever, $kX=\oplus_{\CO}
\CM_{\CO}$
is a decomposition of $kX$ into subrepresentations.
\end{propos}
\proof The action of $\delta_{usu^{-1}}\in k(X)$ on $kX$ in Lemma~5.2
is the projection operator
\[ \delta_{usu^{-1}}\la \phi_{vt}=\delta_{usu^{-1},vtv^{-1}}\phi_{vt}.\]
This is  evident since the action of $k(X)$ is evaluation against
$||\phi_{vt}||$ (or,
explicitly, put $\delta_{usu^{-1}}=\delta_{u(s\la u^{-1})(s\ra u^{-1})}$
into Lemma~5.2.)
We see that $\pi_{[us]}\equiv \delta_{usu^{-1}}\la$ projects $k X$ onto the
subspace
$\CS_{[us]}$ and
\eqn{Ldecomp}{kX=\mathop\oplus_{\chi}\CS_{\chi}
=\mathop\oplus_{\CO}\mathop\oplus_{\chi\in\CO} \CS_{\chi}=
\mathop\oplus_{\CO}\CM_{\CO}   }
as vector spaces. Then the operator
\[
\pi_{\CO} \ =\ \sum_{\chi\in \CO} \pi_{\chi}
\]
commutes with the action of $k(X)\lcross kX$, and is a projection to
$\CM_{\CO}$.
To see that $\pi_{\CO}$ does commute with the action of the algebra, we can
calculate
\[
\pi_{\chi}.(x\tens\delta_y)\ =\ (x\tens\delta_y).\pi_{x^{-1}\tilde\la\chi}\ ,
\]
and note that the operation $x^{-1}\tilde\la$ is a 1-1 correspondence on
the set $\CO$. \endproof

In what follows, we fix an orbit $\CO\subset X/\!\!\sim$ and a
base point $\chi_0$ on it. We denote by
$\Xi$ the isotropy subgroup at $\chi_0$, and $\CS=\CS_{\chi_0}$.

\begin{propos} Let $\CS=\CS_1\oplus \cdots\oplus\CS_n$ be a decomposition  into
irreducibles under the action of $\Xi$. Let $\chi.\CS_i=us\la\CS_i$ when
$\chi=us\tilde\la\chi_0$
(this is independent of the choice of $us$). Then $\CM_i=\oplus_{\chi\in\CO}
\chi.\CS_i\subset \CM_{\CO}$ are
irreducible subrepresentations under $k(X)\lcross kX$. Moreover, $
\CM_{\CO}=\oplus_i \CM_i$ is a decomposition of $\CM_{\CO}$ into irreducibles.
\end{propos}
\proof    Let $x_{\chi}\in X$ be a choice of $us$ such that
$us\tilde\la \chi_0=\chi$. Define
$\chi.\CS_i=x_{\chi}\la\CS_i$.
Now if we take any $x$ so that
$x\tilde \la  \chi_0=\chi$, then $x_{\chi}^{-1}x\in \Xi$, so
\[
x\la \CS_i=x_{\chi}\la(x_{\chi}^{-1}x\la \CS_i)=x_{\chi}\la \CS_i =\chi.\CS_i\
,
\]
as $\CS_i$ is a representation of $\Xi$. Moreover,
$\chi.\CS_i\cup\eta.S_i=\{0\}$ for
$\chi,\eta$ distinct, so $\CM_i$ spanned as shown is a direct sum.

Next we show that $k(X)\lcross kX$ acts on $\CM_i$. Clearly, $x\la \CM_i\subset
\CM_i$ for all $x\in X$ since
$x\la \chi.\CS_i=xx_{\chi}\la \CS_i=\eta.\CS_i$, where
$\eta=x\tilde\la\chi=xx_{\chi}\tilde\la\chi_0$ is another point in $\CO$.
Meanwhile,
The action of the element $\delta_{us}$ in $k(X)$ is either either zero or one
of the
projections associated to (\ref{Ldecomp}), and these are all zero or the
identity on each $\CS_\chi$. Therefore the whole action of $k(X)\lcross kX$
preserves the decompositions of the $\CS_\chi$, with the result that  $\CM_i$
are
subrepresentations and $\CM_{\CO}=\oplus_i \CM_i$.

 We now prove irreducibility. Let $m\in \CM_i$ be nonzero. Then there
is some $\chi$ such that the projection $m_\chi$ to $\chi.\CS_i$ is nonzero,
and then $x_{\chi}^{-1}\la m_\chi$ is a nonzero element of  $\CS_i$.
Since $\CS_i$ is irreducible under $\Xi\subseteq X$, we know that
vectors of the form $\xi x_{\chi}^{-1}\la m_\chi$, for $\xi\in \Xi$,
 span all of $\CS_i$. Since the projection is itself the
action of an element of $k(X)\lcross kX$, we see that $\CS_i$ is contained
in the space spanned by the action of this algebra on $m$.
Moreover by using $x_{\eta}\la \CS_i=\eta.\CS_i$ we see that every
$\eta.\CS_i$ is
contained in the image of $m$ under  $k(X)\lcross kX$. Hence
$\CM_i=\big(k(X)\lcross
kX\big).m$, i.e. $\CM_i$ is irreducible. \endproof

These two propositions give a total decomposition of $kX$ into irreducibles.
In particular, we obtain irreducible subrepresentations for every choice of
orbit and every irreducible subrepresentation of the associated isotropy group.
The converse also holds by similar arguments to those in the preceding
proposition.

\begin{propos} Let $\CM\subset kX$ be an irreducible subrepresentation under
$k(X)\lcross kX$ in Lemma~5.2. Then as a vector space, $\CM$ is of the form
\[ \CM=\mathop\oplus_{\chi\in \CO} \chi.\CM_0\]
for some orbit $\CO$ (with base point $\chi_0$)
 and some irreducible subrepresentation $\CM_0\subset \CS$
under $\Xi$. (Here $\chi.\CM_0=us\la \CM_0$ when $\chi=us\tilde\la \chi_0$.)
\end{propos}
\proof Consider $\CM\subset kX$ and let $\CM_\chi=\pi_{\chi}(\CM)$  for any
$\chi\in
X/\!\! \sim$. Choose a $\chi_0$ so that $\CM_{\chi_0}$ is nonzero,
 write $\CM_0=\CM_{\chi_0}$, and let $\Xi$ be the stabiliser of
$\chi_0$. Now $\CM_0$ must be a representation of $\Xi$. If $\CS_1$ is an
irreducible subrepresentation of $\CM_0$ under $\Xi$, then the previous
proposition shows that
\[
\oplus_{\chi\in \CO} \chi.\CS_1 \subset \CM
\]
is an (irreducible) representation of the algebra,
where $\CO$ is the orbit containing $\chi_0$. Since $\CM$ is irreducible,
this representation
must be equal to $\CM$, and in particular $\CM_0$ is an irreducible
representation
of $\Xi$.
\endproof

It is now a short step to obtain the classification for subrepresentations of
$\ker\eps\subset H$. Note that every Hopf algebra is a direct sum
$k1\oplus\ker\eps$ as
vector spaces, with associated projection $\Pi(h)=h-1\eps(h)$. In our case,
remembering that
\[
\eps(\phi_{vt})=\eps\big( t^{-1}\ra v^{-1}\tens \delta_v\big)=\delta_{v,e}\ ,
\]
we see that
$\ker\eps$ is spanned by the projected basis elements
\[ \bar\phi_{vt}\equiv\Pi(\phi_{vt})=\phi_{vt}-\delta_{v,e}\sum_{u\in
G}\phi_{ue}.\]

\begin{lemma}\label{projact}  The action in Lemma~5.2 and
the projected action
\[ us\la\bar\phi_{vt}=\bar\phi_{us\tilde\la vt},\quad
\delta_{us}\la\bar\phi_{vt}=\delta_{us,vtv^{-1}}\bar\phi_{vt}\]
on $\ker\eps$ are intertwined
by $\Pi:kX\to \ker\eps$.
\end{lemma}
\proof
First we check the $kX$ actions,  using $\Pi$ and then acting in $\ker\eps$
by $us$, to get
\[
\Pi us \la \Pi \phi_{vt}\ =\ \Pi \big(us\la
\phi_{vt}\big)-\delta_{v,e}\Pi(1)\ =\ \Pi \big(us\la \phi_{vt}\big) \ ,
\]
using the facts that $us\la 1=1$ and $\Pi(1)=0$. Now we check the $k(X)$
actions,
using $\Pi$ and then acting in $\ker\eps$ by $\delta_{us}$, to get
\[
\Pi \delta_{us} \la \Pi \phi_{vt}\ =\ \Pi \big(\delta_{us}\la
\phi_{vt}\big)-\delta_{v,e}\Pi(\delta_{us}\la 1)\ =\
\Pi \big(\delta_{us}\la \phi_{vt}\big)-
\delta_{v,e}\Pi(1)\delta_{us,e}\ =\ \Pi \big(\delta_{us}\la \phi_{vt}\big) \ .
\]
\endproof

\begin{theorem}\label{classif} The irreducible quantum tangent spaces $L\subset
\ker\eps$ are all given by the following two cases:

(a)\quad For an orbit $\CO\neq \{[e]\}$ in $X/\!\! \sim$, choose a base
point $\chi_0\in \CO$.
 For each irreducible subrepresentation
$\CM_0\subset \CS$ of $\Xi$ we have an irreducible quantum double
subrepresentation
$\CM=\oplus_{\chi\in\CO}\chi.\CM_0\subset kX$, and an isomorphic quantum double
subrepresentation $L=\Pi
\big(\CM\big)\subset \ker\eps$.

(b)\quad For the orbit $\CO= \{[e]\}$ in $X/\!\! \sim$, choose the base
point $\chi_0=[e]\in \CO$.
For any irreducible subrepresentation
$\CM_0\subset \CS$ of $\Xi$ other than the trivial one (multiples of $1$), we
have an irreducible quantum double subrepresentation
$\CM=\oplus_{\chi\in\CO}\chi.\CM_0\subset kX$, and an isomorphic quantum double
subrepresentation $L=\Pi
\big(\CM\big)\subset \ker\eps$. Here $\CS=kG$ and $\Xi=X$, as in Example~5.4.
\end{theorem}
\proof
In the cases above, $\CM=\oplus_{\chi\in\CO}\chi.\CM_0$ is an irreducible
representation
of the un-projected action. By the previous lemma
$\Pi:\CM\to L$ is a map of representations, where $L\subset \ker\eps$ uses the
projected representation.
The map is onto, and if it is 1-1 then the two representations are
isomorphic, and hence
$L$ is also irreducible.

The only case where the map $\Pi$ is not 1-1 is where $1\in
\CM$. Since $k1$ is a representation and $\CM$ is assumed irreducible, this is
just the nontriviality exclusion stated in the theorem.
We have shown that the cases described lead to irreducible representations of
$\ker\eps$.

Now suppose that we have $L$, an irreducible representation of
$\ker\eps$. Then its inverse image $\Pi^{-1}L\subset H$ is a representation of
$k(X)\lcross kX$, and it contains the subrepresentation $k1$. If $L\neq
\{0\}$, then there
is at least one other irreducible subrepresentation $\CM\subset H$ so that
$k1\oplus \CM\subset
\Pi^{-1}L$. But now $\CM$ must be of the form described earlier, and by
irreducibility
$\Pi \CM=L$.
\endproof

We are now left with the problem of finding the irreducible
subrepresentations of
$\CS$ under the group $\Xi$. We can decompose the representation
$\CS$ into a sum of irreducibles $\CS=\CS_1\oplus\dots\oplus\CS_r$, so in
this case there would be at least
$r$ irreducible subrepresentations. However if there are any equivalent
representations in this list,
there are many more possible irreducible subrepresentations. This is a standard
situation in representation theory and we briefly recall its resolution, as
follows.

\begin{lemma}
Let $\CN$ be an irreducible representation of the group $\Xi$, and let
 $\CM\subset \CN\oplus U$ be an irreducible representation of $\Xi$, for
another representation
$U$. Suppose
 that there is a vector in $\CM$ with a nonzero $\CN$ component. Then there is
a $\Xi$-map $T:\CN\to U$
so that $\CM=\big\{w\oplus Tw |\  w\in \CN\big\}$.
\end{lemma}
\proof
First we show that any $w\in \CN$ occurs as the first coordinate of a vector
in $\CM$. The projection
$\pi:\CN\oplus U\to \CN$ is a $\Xi$-map, so $\pi \CM\subset \CN$ is a
subrepresentation of $\CN$. As $\CN$ is
irreducible, we see that $\pi \CM = \CN$.

Now we
show that for any $w\in \CN$ there is at most one $u\in U$ for which $w\oplus
u\in \CM$. If we had $w\oplus u\in \CM$
and $w\oplus u'\in \CM$, then on subtraction we would also have
$0\oplus(u-u')\in \CM$. But now using
$0\oplus(u-u')$ as a cyclic vector we could construct a representation
contained in $\CM$ which
consisted purely of vectors with zero $\CN$ component. This would contradict
the irreducibility of
$\CM$ unless $u=u'$.

Now simply label this unique $u$ as $T(w)$. As $\xi(w\oplus Tw)=\xi w\oplus
\xi Tw\in \CM$ for
$\xi \in \Xi$, we must have $T(\xi w)=\xi T(w)$.
\endproof

\begin{propos} Suppose that $\CM$, an  irreducible representation of the
group $\Xi$, is
contained in $\CN_1^{\oplus n_1}\oplus\dots\oplus \CN_k^{\oplus n_k}$, where
the
$\CN_i$ are inequivalent irreducible representations, and $n_i$ gives the
corresponding multiplicities
in the sum. Then there is an index $j$, and
$(\lambda_1,\dots,\lambda_{n_j})\in k^{n_j}$
(not all zero)
so that
\[
\CM\ =\ \Big\{ 0\oplus\dots\oplus 0\oplus (\lambda_1 w\oplus\dots \oplus
\lambda_{n_j} w)\oplus
0\oplus\dots\oplus 0 \in \CN_1^{\oplus n_1}\oplus\dots\oplus \CN_k^{\oplus n_k}
|\
w\in \CN_j\Big\}\ .
\]
\end{propos}
\proof
A nonzero vector in $\CM$ must have a nonzero component in at least one of the
irreducible summands.
To reduce confusion, we shall assume that the first $\CN_1$ summand has a
such a
nonzero component. Then we can write $\CM\subset \CN_1\oplus U$, where $U$ is
the sum of the
rest of the $\CN_i$s. By the previous lemma, there is a $\Xi$-map $T:\CN_1\to
U$ so that
$\CM=\{w\oplus Tw|\ w\in \CN_1\}$.

Now we can follow the map $T$ with a projection to
one of the $\CN_i$ summands in $U$ and examine the resulting map.
If we project to another of the $\CN_1$ summands, the result is a $\Xi$-map
from an
irreducible representation to itself, and by Schur's lemma this must be
scalar multiplication.
If we project to a $\CN_i$ summand for $i\neq 1$, the result is a $\Xi$-map
from an
irreducible representation to an inequivalent
irreducible representation, and by Schur's lemma this must be zero.
\endproof

Thus the irreducible subrepresentations of $\CS$ in Theorem~5.9 are themselves
classified as follows. For a given orbit $\CO$ with base point $\chi_0$,
decompose $\CS$ into irreducible representations under the action of $\Xi$.
Write this as
$\CS=\CM_1^{\oplus n_1}\oplus\dots\oplus \CM_k^{\oplus n_k}$, where the
$\CM_i$ are inequivalent irreducible representations, and $n_i$ gives the
corresponding multiplicities
in the sum. The data is then $(\CO,\CM_0,\lambda)$, where $\CM_0$ is an
irreducible representation
 of $\Xi$ occuring in $\CS$
with multiplicity $n>0$, and $\lambda\in \Bbb P(k^n)$. (Note that we take
the projective space as only
 the ratios of the scalars in the last proposition are required to specify
the subspace.)

Bicrossproducts interpolate between group algebras and group function algebras.
As a check, we recover the seemingly quite disparate results
known separately for these two cases.

\begin{example} Suppose that $G$ is trivial. Then $H=kM$. In this case $\sim$
is the same as equality and $X/\!\! \sim=X$. In this case the equivalence
classes
are singletons corresponding to points in $M$ and $s\tilde\la t=sts^{-1}$ is
conjugation. Hence the orbits $\CO$ are conjugacy classes in $M$. The action of
the isotropy group is trivial and hence this is the only data. We recover the
result \cite{BMD:non}\cite{Ma:cla} that the irreducible quantum tangent spaces
are the
projected spans of the conjugacy classes.
\end{example}

\begin{example} Suppose that $M$ is trivial. Then $H=k(G)$. In this case $X$ is
one entire equivalence class and $X/\!\! \sim=\{[e]\}$. Then there is only
one orbit
$\CO={[e]}$ and Theorem~\ref{classif} reduces to the classification of quantum
tangent spaces in $\ker\eps\subset k(G)$ (so $\CO=kG$) which have been
classified by the second author in \cite{Ma:cla} as irreducible subspaces under
the left regular representation induced by $u\tilde\la v=uv$.
\end{example}

We can also keep $M,G$ non trivial but let $\la$ or $\ra$ be trivial in the
matched pair. In these cases $X$ is a semidirect product.
 Note that the two cases are quite different. In one, $H$ is a tensor
product algebra. In the other it is a cross product $kM\rcross k(G)$ where $M$
acts by group automorphisms of $G$.

\begin{example} Suppose that the action $\la$ is trivial.
If we start with an orbit $\CO$ containing $[s]$, then
$[s]=sG$ and $\Xi_{[s]}\ =\ \big\{ ut|\ tst^{-1}=s\ra u\big\}$. The canonical
projection $ut\mapsto u$ us a group homomorphism $\Xi\to G$ and, identifying
$\CS=kG$, the action
of $\Xi$ is the pull back along this of the left regular representation.
\end{example}

\begin{example} Suppose that the action $\ra$ is trivial.
If we start with an orbit $\CO$ containing $[s]$, then
$[s]\ =\ \big\{ sv|\ sv=vs\big\}$ may be identified with $G_s$, the $\cup G$ of
the centraliser of $s\in X$.
The isotropy group is $\Xi_{[s]}\ =\ \big\{ ut|\ su=us\ \ {\rm and}\ \
st=ts\big\}=G_s\times M_s$, where $M_s$ is the centraliser of $s\in M$.
Identifying $\CS=k G_s$, the action of $\Xi$ is with $M_s$ acting by $\la$ and
$G_s$ acting by the left regular representation.
\end{example}

\begin{example} We completely solve the problem by hand for
Example~2.6 where $X=S_3$.  We just
decompose $k(X)$ into irreducible subspaces $L$ under (\ref{projact}) and
project to $\ker\eps$.
More precisely, we follow the above theory in this special case and then
convince ourselves by direct means that it is the right answer.
\end{example}

{\bf The orbit decomposition}\quad  First we find the allowed orbits:

Orbit 1 : $\CO=\{[e]\}$, choosing base point $[e]$. Then
$[e]=\big\{e,(123),(321)\big\}$ and $\Xi=S_3$. The action of
$\Xi$ on $[e]$ is given by the formula $ut\tilde\la v=utvt^{-1}$. We have the
eigenspaces of the $G$ action $\CM_0=\big< e+(123)+(321)\big>$,
$\CM_1=\big< e+\omega(123)+\omega^2(321)\big>$, and
$\CM_2=\big< e+\omega^2(123)+\omega(321)\big>$. Now $\CM_0$ forms a
1-dimensional representation of $\Xi$,
but this is annihilated by $\Pi$. The action of $(12)\in \Xi$ is to swap
$\CM_1$ and $\CM_2$, so we get
a 2 dimensional irreducible representation $\CM_1\oplus \CM_2$ of $\Xi$, giving
a
2 dimensional irreducible subrepresentation in $\ker\eps$.

Orbit 2 : $\CO=\{[(12)],[(13)],[(23)]\}$, choosing base point $[(12)]$.
Then $[(12)]=\big\{(12)\big\}$
and $\Xi=M$. The action of $\Xi$ on $[(12)]$ is the trivial 1 dimensional
representation, giving a 3 dimensional irreducible representation in
$\ker\eps$.

{\bf The direct approach}\quad  The space $\ker\eps$ is spanned by
the vectors $\big\{\bar\phi_{x}|\ x\in S_3\big\}$, where
$\bar\phi_x=\Pi\phi_x$,
and there is the linear relation
$\bar\phi_{e}+\bar\phi_{(123)}+\bar\phi_{(321)}=0$. We use the relation to
rewrite $\bar\phi_{(321)}=-\bar\phi_{e}-\bar\phi_{(123)}$, giving
a basis consisting of 5 elements. This is then split into two parts by the
action of
$X$:

(1)\quad The space spanned by the elements
$\big\{\bar\phi_{e},\bar\phi_{(123)}\big\}$.
This has the $X$ action $us\la \bar\phi_{e}=\bar\phi_{u}$ and
$us\la \bar\phi_{(123)}=\bar\phi_{us(123)s^{-1}}$, where we remember to rewrite
$\bar\phi_{(321)}=-\bar\phi_{e}-\bar\phi_{(123)}$. This gives a 2 dimensional
irreducible representation.

(2)\quad The space spanned by the elements
$\big\{\bar\phi_{(12)},\bar\phi_{(13)},\bar\phi_{(23)}\big\}$.
This has the $X$ action $us\la \bar\phi_{v(12)}=\bar\phi_{usvs^{-1}(12)}$
for all $v\in G$. This gives a 3 dimensional
irreducible representation.

\begin{example} We apply the theory above to Example~3.4 where $X=C_6\dcross
C_6$ and
spell out the final result. This is rather complicated to do directly
(hence justifying our methods). It is one of the simplest
examples\cite{BGM:fin}
of a true bicrossproduct Hopf algebra.
\end{example}

First we identify the possible values of $\|.\|$, and the orbits of these
values under $X=S_3
\times S_3$.  Since $\|vt\|=vtv^{-1}$, the possible values and the orbits are
simply given by looking at the conjugacy classes of the elements $t\in M$
in $X$. These are:

{\bf (Orbit 0)}\quad $t=0_M$ gives the conjugacy class consisting only of
the identity.
We choose $[0_M]$ as the
base point for this orbit.  Then
$\Xi_{[0]}=X$, and the equivalence class is $[0_M]=G$ (as previously noted).
The action of $\Xi_{[0]}=X$ on $[0_M]=G$ is given by the formula
$ut\tilde\la v=u(t\la v)$.

Let us now look at the decomposition of $\CS=kG$ into irreducibles under
the action of
$\Xi_{[0]}=X$. The element $1_G\in \Xi$ acts on any irreducible $\CM\subset
\CS$, and
its action diagonalises, that is $\CM$ is a sum of $\CM_r$ ($r\in \Bbb Z_6$),
where each $\CM_r$ is zero or
$\big<f_r\big>_k$, where
\[
f_r\ =\
0_G+\omega^{r}1_G+\omega^{2r}2_G+\omega^{3r}3_G+\omega^{4r}4_G+\omega^{5r}5_
G\ ,
\]
for $\omega$ a
primitive $6$th root of unity. The action of $1_M$ is to send $f_r$ to
$f_{-r}$, so we get
the 4 irreducible representations $\CS_1=\big<f_0\big>$,
$\CS_2=\big<f_1,f_5\big>$, $\CS_3=\big<f_2,f_4\big>$, and
$\CS_4=\big<f_3\big>$.
The two 1-dimensional representations $\CS_1$ and $\CS_4$ are not equivalent,
as
$\CS_1$ is the trivial representation and $\CS_4$ is not trivial. To show that
$\CS_2$ and $\CS_3$ are not equivalent we use the trace of $1_G$ on the
representations,
which is $\omega+\omega^{-1}$ on $\CS_2$ and $\omega^2+\omega^{-2}$ on $\CS_3$.

There are four inequivalent irreducible representations for this orbit, but
one is annihilated by $\Pi$, leaving three irreducible representations
of $\ker\eps$ on application of $\Pi$.

{\bf (Orbit 1)}\quad $t=1_M$ and $t=5_M$ give the conjugacy class consisting
of elements of the form (any 2-cycle in 1,2,3)(any 3-cycle in 4,5,6). We
choose $[1_M]$ as the
base point for this orbit. Then $\Xi_{[1]}=M$, and the equivalence class is
$[1_M]=\big\{ 1_M v|\ 1_M\la v=v\big\}$. Since $1_M$ acts on $G$ by the
permutation
$(1_G,5_G)(2_G,4_G)$, we find that $[1_M]=\big\{ 1_M 0_G,1_M 3_G\big\}$.
The action of $\Xi_{[1]}=M$ on $[1_M]$ is given by the formula $t\la 1_M
v=1_M (t\la v)$,
 which is the trivial action since both $0_G$ and $3_G$ are fixed by the left
action of $M$.

The decomposition of $\CS=k\big\{ 1_M 0_G,1_M 3_G\big\}$ into irreducibles
under the action of
$\Xi_{[1]}=M$ gives two copies of the trivial one-dimensional representation.

{\bf (Orbit 2)}\quad $t=2_M$ and $t=4_M$ give the conjugacy class consisting
of elements of the form (any 3-cycle in 4,5,6). We choose $[2_M]$ as the
base point for this orbit. Then $\Xi_{[2]}=S_3\times C_3$, where $C_3$ is
the group
of permutations of $\{4,5,6\}$ consisting of the identity and the two 3-cycles.
In terms of the factorisation, $\Xi_{[2]}$ consists of all elements of the
form $ut$
for $u\in\big\{0_G,2_G,4_G\big\}$ and $t\in M$.
The equivalence class is $[2_M]=\big\{ 2_M v|\ 2_M\la v=v\big\}$.
Since $2_M$ has the trivial action on $G$, we find that
$[2_M]=\big\{ 2_M 0_G,2_M 1_G,2_M 2_G,2_M 3_G,
2_M 4_G,2_M 5_G\big\}$. Under the $\Xi_{[2]}$ action $ut\tilde sv=su(t\la
v)$, $[2_M]$
splits into two orbits,
$\big\{ 2_M 0_G,2_M 2_G,
2_M 4_G\big\}$ and
$\big\{ 2_M 1_G,2_M 3_G,2_M 5_G\big\}$.

First we decompose $k\big\{ 2_M 0_G,2_M 2_G,2_M 4_G\big\}$ into
 irreducibles under the action of
$\Xi_{[2]}$. The action of $2_G$ on this vector space
diagonalises, that is $k\big\{ 2_M 0_G,2_M 2_G,2_M 4_G\big\}=\CM_0\oplus
\CM_1\oplus \CM_2$,
where $\CM_0=\big<2_M 0_G+2_M 2_G+2_M 4_G\big>_k$,
$\CM_1=\big<2_M 0_G+\omega 2_M 2_G+\omega^2 2_M 4_G\big>_k$, and
$\CM_2=\big<2_M 0_G+\omega^2 2_M 2_G+\omega 2_M 4_G\big>_k$ ($\omega$ being a
primitive 3rd root of unity).
 The effect of $1_M$ on these eigenspaces is to swap $\CM_1$
and $\CM_2$. The decomposition into irreducibles gives $\CS_1=\CM_0$
(trivial 1-dimensional) and $\CS_2=\CM_1\oplus \CM_2$.

Next we decompose $k\big\{ 2_M 1_G,2_M 3_G,2_M 5_G\big\}$ into
 irreducibles under the action of
$\Xi_{[2]}$.
The action of $2_G$ on this vector space
diagonalises, that is $k\big\{ 2_M 1_G,2_M 3_G,2_M 5_G\big\}=\CN_0\oplus
\CN_1\oplus \CN_2$,
where $\CN_0=\big<2_M 1_G+2_M 3_G+2_M 5_G\big>_k$,
$\CN_1=\big<2_M 1_G+\omega 2_M 3_G+\omega^2 2_M 5_G\big>_k$, and
$\CN_2=\big<2_M 1_G+\omega^2 2_M 3_G+\omega 2_M 5_G\big>_k$.
 The effect of $1_M$ on these eigenspaces is to swap $\CN_1$
and $\CN_2$. The decomposition into irreducibles gives $\CS_3=\CN_0$
(trivial 1-dimensional) and $\CS_4=\CN_1\oplus \CN_2$.

In fact the two 2-dimensional representations
$\CS_2$ and $\CS_4$ are isomorphic, using the map
$2_M 0_G+\omega 2_M 2_G+\omega^2 2_M 4_G \mapsto 2_M 1_G+\omega 2_M
3_G+\omega^2 2_M 5_G$
and
$2_M 0_G+\omega^2 2_M 2_G+\omega 2_M 4_G \mapsto \omega^2(
2_M 1_G+\omega^2 2_M 3_G+\omega 2_M 5_G)$.

{\bf (Orbit 3)}\quad $t=3_M$ gives the conjugacy class consisting
of elements of the form (any 2-cycle in 1,2,3). We choose $[3_M]$ as the
base point for this orbit. Then $\Xi_{[3]}=C_2\times S_3$, where $C_2$ is
the group
of permutations of $\{1,2,3\}$ consisting of the identity and $(1,2)$.
In terms of the factorisation, $\Xi_{[3]}$ consists of all elements of the
form $ut$
for $u\in\big\{0_G,3_G\big\}$ and $t\in M$.
The equivalence class is $[3_M]=\big\{ 3_M v|\ 3_M\la v=v\big\}$.
Since $3_M$ acts on $G$ by the permutation
$(1_G,5_G)(2_G,4_G)$, we find that $[3_M]=\big\{ 3_M 0_G,3_M 3_G\big\}$.
The $\Xi_{[3]}$ action on $[3_M]$ is given by $ut\tilde\la sv=su(t\la v)=suv$.

Now decompose $\CS=k\big\{ 3_M 0_G,3_M 3_G\big\}$ into
 irreducibles under the action of
$\Xi_{[3]}$.  The action of $3_G$ on this vector space
diagonalises, that is $\CS=\CM_1\oplus \CM_2$,
where $\CM_1=\big<3_M 0_G+3_M 3_G\big>_k$ and $\CM_1=\big<3_M 0_G-3_M
3_G\big>_k$.
The decomposition into irreducibles gives $\CS_1=\CM_1$ and $\CS_2=\CM_2$,
which are not equivalent.

This completes our classification of the bicovariant quantum tangent spaces on
bicrossproducts. It remains to dualise the results to obtain the corresponding
first order differential
calculi, as outlined at the start of the section. Explicitly, the dual of the
inclusion
$i_L:L\to \ker\eps_H$ is a surjection $i_L^*:(\ker\eps_H)^*\to V$,
where $V=L^*$. On the other hand, the inclusion $j:\ker\eps_H\to H$ dualises to
a map
$j^*:A\to(\ker\eps_H)^*$ where $A=H^*$. Since $\ker j^*=(\image
j)^{\perp}=(\ker\eps_H)^{\perp}=k 1_{A}$,
the restriction $j^*|_{\ker\eps_{A}}:\ker\eps_{A}\to (\ker\eps_H)^*$ is
an isomorphism.
Putting these together, we get a quotient map $\pi_V:\ker\eps_{A}\to V$.
Recall also that we can describe a
representation $L$ of
$D(H)$ as a left $H$ and right $H^*$-module,  with actions obeying the
compatibility
condition
\[
 h\la (x\ra a)=\sum( (h\o\ra a\o)\la x) \ra(h\t\la a\t)\ ,
\]
for all $x\in L$, which can be further computed in terms of the mutual
coadjoint actions. (We
freely identify a left $H^{*\rm op}$ module as a right $H^*$-module.) For a
given $a\in H^*$, the action $\ra a:L\to L$ dualises
to
$(\ra a)^*:V\to V$ where $V=L^*$. Similarly for the operators $h\la$. We obtain
in this way a left action of $A=H^*$ and a right action of $H$ on $V$, by
$a\la v=(\ra a)^*(v)$
and $v\ra h=(h\la)^*(v)$, which make $V$ into a $D(A)$ representation. This is
a general observation about the dualisation of quantum double modules. Combined
with the projection $\pi_V$ and $\extd$ in (\ref{extd}), we obtain the
corresponding first order differential calculus $\Omega^1=V\tens A$.

To fully specify the resulting exterior differential $\extd$ it is equivalent
and more
convenient to give its evaluations $\del_v=(\<v,\ \>\tens\id)\circ \extd$
against all $v\in L$. These operators $\del_v:A\to A$ are the {\em braided
vector fields} associated to
elements $v\in L$. The term is justified because they obey a braided version of
the Leibniz rule\cite{Ma:cla}. Since we obtain $L$ from the map $\Pi:\CM\to L$,
we specify equivalently
the operators
\[ \del_m=(\<m,\ \>\tens\id)(\Pi^*\tens\id)\circ\extd,\quad\forall m\in \CM.\]

\begin{example} The $\del$ operators in a few cases from Example~5.17.
\end{example}

We take orbit 0 in Example~5.17, and write
\[
f_r\ =\ 0_M\tens\delta_{0_G}+\omega^{r} 0_M\tens\delta_{1_G}
+\omega^{2r} 0_M\tens\delta_{2_G}+\omega^{3r} 0_M\tens\delta_{3_G}
+\omega^{4r} 0_M\tens\delta_{4_G}+\omega^{5r} 0_M\tens\delta_{5_G}\ ,
\]
where $\omega$ is a primitive 6th root of unity. From its definition,
\[\del_{f_r}(\delta_t\tens v)\ =\ \sum_{ab=t} \big<\delta_a\tens b\la
v,f_r\big>\ \delta_b\tens v\ -\
\sum_{a} \big<\delta_a\tens 0_G,f_r\big>\ \delta_t\tens v\ ,
\]
After a little calculation this reduces to
\[\del_{f_r}(\delta_t\tens v)\ =\ (\omega^{r(t\la v)}-1)\
\delta_t\tens v\ .
\]
We use the allowed spaces for $\CM$ for this orbit, which are
$\big<f_3\big>_k$, $\big<f_1,f_5\big>_k$, and $\big<f_2,f_4\big>_k$.

Next we take a case from orbit 2, the irreducible representation
given by $\CM_0=\big<g_1,g_2\big>_k$, where
\[
g_r\ =\ \phi_{2_M0_G}+\omega^r \phi_{2_M2_G}+\omega^{2r}\phi_{2_M4_G}\ =\
4_M\tens\delta_{0_G}+\omega^{r}4_M\tens\delta_{2_G}+
\omega^{2r}4_M\tens\delta_{4_G}
\ ,
\]
and $\omega$ is a primitive 3rd root of unity.
In this case the orbit consists of more than one point, in fact $\CO=\big\{
[2_M],[4_M]\big\}$.
We choose $x_{[4_M]}\in X$ so that $x_{[4_M]}\tilde\la [2_M]=[4_M]$, for
example
$x_{[4_M]}=1_G$. Now we can add $1_G \tilde\la g_1$ and $1_G \tilde\la g_2$
to $g_1$ and $g_2$
to get a basis of the 4-dimensional representation specified by $\CO$ and
$\CM_0$, that is
$\CM=\big<g_1,g_1,1_G \tilde\la g_1,1_G \tilde\la g_2\big>_k$, where
\[
1_G \tilde\la g_r\ =\
\phi_{1_G2_M}+\omega^r \phi_{3_G2_M}+\omega^{2r}\phi_{5_G2_M}\ =\
2_M\tens\delta_{1_G}+\omega^{r}2_M\tens\delta_{3_G}+\omega^{2r}2_M\tens
\delta_{5_G}\ .
\]
Then we may calculate the evaluations
\[\del_{g_r}(\delta_t\tens v)\ =\ \Big( \delta_{(t-4)\la v,0}+\omega^r
\delta_{(t-4)\la v,2}+\omega^{2r} \delta_{(t-4)\la v,4}\Big)\
\delta_{t-4}\tens v\ -\ \delta_{t}\tens v\ ,
\]
\[\del_{1_G \tilde\la  g_r}(\delta_t\tens v)\ =\ \Big( \delta_{(t-2)\la
v,1}+\omega^r
\delta_{(t-2)\la v,3}+\omega^{2r} \delta_{(t-2)\la v,5}\Big)\
\delta_{t-2}\tens v\ .
\]


\end{document}